\documentclass[11pt]{article}

\usepackage{latexsym}
\usepackage{bbm}
\usepackage{amsmath}
\usepackage{amsfonts}
\usepackage{amssymb}
\usepackage{multirow}
\usepackage{latexsym}
\usepackage[dvips]{graphicx}
\usepackage{overpic}
\usepackage{graphicx}
\usepackage{relsize}
\usepackage{leftidx}
\allowdisplaybreaks
\newtheorem{theorem}{\bf Theorem}[section]
\newtheorem{lemma}[theorem]{\bf Lemma}

\newtheorem{defn}{\bf Definition}[section]

\newtheorem{remark}{{\bf Remark}}[section]
\newenvironment{proof}{\noindent{\em Proof.}}{\quad \hfill$\Box$\vspace{2ex}}

\def \bZ {\Bbb Z}

\def \bR {\Bbb R}

\def \and {\, \mb
ox{\rm and}\, }

\def \supp {\,{\rm supp}\,}

\def \Re {\,{\rm Re}\,}

\def \min{\,{\rm min}}
\def \max{\,{\rm max}}

\def \l {\left}
\def \r {\right}
\makeatletter

\newcommand{\Rmnum}[1]{\expandafter\@slowromancap\romannumeral #1@}

\makeatother
\begin{document}

\title {\bf Damping estimates for oscillatory integral operators with real-analytic phases and its applications}

\author{Zuoshunhua Shi \thanks{School of Mathematics and Statistics, Central South University, Changsha, People's Republic of China. E-mail address:
        {\it shizsh@163.com}.}\,\,\,,\quad
       Shaozhen Xu\thanks{School of Mathematical Sciences, University of
        Chinese Academy of Sciences, Beijing 100049, People's Republic of China. \& Department of Mathematics, University of Illinois at Urbana-Champaign, Urbana, IL, 61801, USA. E-mail
        address: {\it xushaozhen14b@mails.ucas.ac.cn}.}
         \quad and \quad Dunyan Yan\thanks{School of Mathematical Sciences, University of Chinese Academy of Sciences, Beijing 100190, P. R. China. E-mail address: {\it ydunyan@ucas.ac.cn}.} }

\date{}
\maketitle{}

\begin{abstract}
In this paper, we investigate sharp damping estimates for a class of one dimensional oscillatory integral operators with real-analytic phases. By establishing endpoint estimates for suitably damped oscillatory integral operators, we are able to give a new proof of the sharp $L^p$ estimates which have been proved by Xiao in Endpoint estimates for one-dimensional oscillatory integral operators, \emph{Advances in Mathematics}, \textbf{316}, 255-291 (2017). The damping estimates obtained in this paper are of independent interest.
\end{abstract}


\section{Introduction}
In this paper, we consider one-dimensional oscillatory integral operators of the form
\begin{equation}\label{section 1 def Tlambda}
T_{\lambda}f(x)=\int_{-\infty}^{\infty} e^{i\lambda S(x,y)} \varphi(x,y) f(y)dy,
\end{equation}
where $\lambda\in\bR$ is a parameter, $S$ is a real-analytic function near the origin in $\bR^2$, and $\varphi$ is a smooth cut-off function supported in a small neighborhood of the origin. Our main goal is to determine the optimal decay of the operator norm of $T_{\lambda}$ on $L^p$. It is easy to see that $T_{\lambda}$ has no power decay property with respect to $\lambda$, if the phase function $S$ is degenerate in the sense that it can be written as
\begin{equation}\label{section 1 degenerate phase}
S(x,y)=P(x)+Q(y)
\end{equation}
for two polynomials $P$ and $Q$. A natural question arises whether there is a power decay estimate if the phase $S$ is not of the form (\ref{section 1 degenerate phase}). The answer is affirmative except the special cases $p=1$ and $p=\infty$. We now review some well-known results in this direction. If $S$ is nondegenerate in the sense that its Hessian $S_{xy}''$ does not vanish in the support of $\varphi$, H\"{o}rmander \cite{hormander2}
obtained the maximal $L^2$ decay $|\lambda|^{-1/2}$. For general real-analytic phases, Phong-Stein \cite{PS1997} proved the remarkable result that the maximal $L^2$ decay for $T_{\lambda}$ is explicitly determined by the Newton polyhedron of $S$; see also Greenblatt \cite{greenblatt1} for a new proof of this result. The Phong-Stein theorem was extended to the case of almost all smooth phase functions by Rychkov \cite{rychkov}. Its full generalization to general smooth functions was established by Greenblatt \cite{greenblatt2}.

Recently, Xiao \cite{Xiao2017} has proved the sharp $L^p$ decay estimates for $T_{\lambda}$. For convenience, we formulate this result in the following theorem. Another equivalent formulation in terms of the Newton polyhedron of $S$ will be given in Section 3.
\begin{theorem}\label{section 1 Main theorem}
{\rm ( \textbf{Xiao \cite{Xiao2017}} ) } Assume $S$ is a real-analytic function near the origin in the plane. Let $T_{\lambda}$ be the oscillatory integral operator defined by {\rm (\ref{section 1 def Tlambda})}. If there exist two positive integers $k,l$ such that $\partial_x^k\partial_y^lS(0,0)\neq 0$, then there is a constant $C=C(S,\varphi)$ such that
\begin{equation}\label{section 1 main estimate}
\|T_{\lambda}f\|_{L^p}
\leq
C|\lambda|^{-\frac{1}{k+l}}\|f\|_{L^p},~~~p=\frac{k+l}{k},
\end{equation}
provided that $\varphi$ is supported in a sufficiently small neighborhood of the origin.
\end{theorem}

By interpolation, one can verify that the estimate (\ref{section 1 main estimate}) is sharp only if $(k,l)$ lies on the boundary of the Newton polyhedron $N(S)$. Moreover, it suffices to show that the above $L^p$ estimates are true for extreme points $(k,l)\in N(S)$. For further results, one can see Greenleaf-Seeger \cite{greenleafseeger1}, Yang \cite{yangchanwoo, yangchanwoo2} and Shi-Yan \cite{ShiYan} for earlier work. For linear and multi-linear estimates, we refer the reader to Carbery-Christ-Wright \cite{carbery}, Carbery-Wright \cite{CaW}, Phong-Stein-Sturm \cite{PSS2001}, Christ-Li-Tao-Thiele \cite{Christ-Li-Tao-Thiele} and Gressman-Xiao \cite{gressmanxiao}. Some work on higher dimensional oscillatory integral operators can be found in Tang \cite{tangwan}, Greenleaf-Pramanik-Tang \cite{GPT} and Xu-Yan \cite{xuyan}. Other results concerning regularity of Radon transforms associated with $S$, we refer the reader to Greenleaf-Seeger \cite{greenleafseeger2} and Seeger \cite{seeger1,seeger2}.

Our proof of Theorem \ref{section 1 Main theorem} was inspired by Xiao \cite{Xiao2017} and the recent work in \cite{Shi} of the first author. It was found in \cite {Shi} that some uniform damping estimates for $T_{\lambda}$ imply sharp uniform $L^p$ decay estimates. In this respect, damping estimates for oscillatory integral operators are of independent interest.

Let us discuss some basic facts related to damped oscillatory integral operators of the form
\begin{equation}
W_zf(x)=\int_{-\infty}^{\infty} e^{i\lambda S(x,y)} |D(x,y)|^z \varphi(x,y) f(y) dy,
\end{equation}
where the damping factor $D$ is determined by the Hessian $S_{xy}''$, and the exponent $z\in\mathbb{C}$ lies in some bounded strip $a\leq \Re(z) \leq b$. On the one hand, we shall give sharp $L^2$ decay estimates for $W_z$ when $z$ has real part $\Re(z)=b$. On the other hand, for $z$ with $\Re(z)=a$, we have to establish $L^1\rightarrow L^{1,\infty}$, $H^1_E\rightarrow L^1$ and $L^1\rightarrow L^{1}$ estimates. In the latter case, the argument is somewhat different depending on whether $\{(x,y):D(x,y)=0\}$ contains at least two curves. Assume $\{D=0\}$ includes only one curve. If this curve is tangent to some straight line at the origin, then we shall define a class of variants of Hardy spaces $H^1$ to replace $L^1\rightarrow L^1$ by $H_E^1\rightarrow L^1$. This idea was explored by Phong-Stein \cite{PS1986}, Greenleaf-Seeger \cite{greenleafseeger1} and Pan \cite{pan}. An interpolation of these two sharp decay estimates will lead to our desired sharp $L^p$ estimates.

This paper will be organized as follows. Some useful lemmas will be presented in Section 2. We shall establish sharp $L^2$ decay estimates for a class of damped oscillatory integral operators in Section 3. In Section 4, we shall prove some endpoint estimates of form $L^1\rightarrow L^{1,\infty}$, $H_E^1\rightarrow L^1$ and $L^1\rightarrow L^1$. A new proof of Theorem \ref{section 1 Main theorem} will be given in Section 5. Throughout this paper, we use the notation $a\lesssim b$ to mean $a\leq Cb$ for an appropriate constant $C>0$, and $a\approx b$ to mean $a\lesssim b$ and $b\lesssim a$.

\section{Preliminaries}
In this section, we shall first give the concept of horizontally (vertically) convex domains which were introduced by Phong-Stein-Sturm \cite{PSS2001}; see also Carbery-Wright \cite{CaW} for some related remarks. This notion of convexity turns out to be important in this paper. In fact, the operator van der Corput lemma will be established for oscillatory integral operators supported on horizontally and vertically convex domains. Finally, we shall give a variant of Stein-Weiss interpolation with change of measures. Some related results will be also included in this section.
\begin{defn}
Let $\Omega$ be a domain in $\bR^2$. If $(x,z), (y,z)\in \Omega$ imply
$(\theta x+(1-\theta)y, z)\in \Omega$ for all $0\leq\theta\leq 1$, then $\Omega$ is said to be horizontally convex. Similarly, a domain $\Omega$ is said to be vertically convex
if $(x,y), (x,z)\in \Omega$ imply $(x,\theta y+(1-\theta)z)\in \Omega$ for all $0\leq\theta\leq 1$.
\end{defn}
It is clear that a convex domain is both horizontally convex and vertically convex. In subsequent sections, we shall use frequently a special horizontally and vertically convex domain which is known as the curved trapezoid; see Phong-Stein \cite{PS1997,PS1998} for earlier work where the curved trapezoid was used as an important domain.
\begin{defn}
Suppose $g$ and $h$ are two monotone functions on $[a,b]$ such that $g(x)\leq h(x)$ for every $x\in[a,b]$. Then the following domain
\[
\Omega=\left\{(x,y):a\leq x\leq b, g(x)\leq y\leq h(x)\right\}
\]
is said to be a curved-trapezoid.
\end{defn}

Now we recall the concept of polynomial type functions; see Phong-Stein \cite{PS1997,PS1998}.

\begin{defn}\label{sec2 def poly type fun}
Assume $F\in C^{N}$ is defined on a bounded interval $J$. Then $F$ is said to be of polynomial type with order $N$ if there exists a constant $C_F>0$ such that
\begin{equation*}
\sup_{x\in J}|F^{(N)}(x)|
\leq
C_F \inf_{x\in J}|F^{(N)}(x)|.
\end{equation*}
We use $ord(F)$ to denote the order of $F$.
\end{defn}

Polynomial type functions have the following useful property; see Phong-Stein \cite{PS1997,PS1998}.
\begin{lemma}\label{sec2 esti poly tpype func}
Let $J$ be a bounded interval. If $F\in C^{N}(J)$ be a polynomial type function with order $N$, then there is a constant $C=C(N,C_F)$ such that
\[
\sum_{k=0}^{N}\l|I^{*}\r|^k\sup_{x\in I^{*}\cap J}\left|F^{(k)}(x)\right|\leq C\sup_{x\in I}|F(x)|
\]
is true for all subintervals $I\subseteq J$, where for each interval $I$ we denote by $I^{*}$ its double, i.e., the interval with the same center as $I$ but with twice its length.
\end{lemma}
\begin{remark}\label{Rmk polynomial type func}
The following basic fact turns out to be useful in subsequent sections.

$\bullet $ Assume $F$ is a real-valued polynomial type function on some bounded interval $J$. If there exist two constants $\mu>0,C>0$ such that $\mu\leq |F(x)|\leq C\mu$ on $J$, then $|F|^z$ also satisfies the estimate in Lemma
{\rm\ref{sec2 esti poly tpype func}} for all complex numbers $z$ with a constant $C=C(ord(F),C_F,z)$.
\end{remark}

With the above concepts, we are able to give the operator van der Corput lemma. The following notations will be frequently used in this
section. Assume $\Omega$ is a horizontally and vertically convex domain. We define the following notations associated with $\Omega$.

\begin{itemize}
\item $\delta_{\Omega, h}(x)$: the length of the cross section $\{y\mid (x,y)\in \Omega\}$ with the subscript $h$ indicating that
$\delta_{\Omega, h}$ is a function of the horizontal component.
\item $\delta_{\Omega,v}(y)$: the length of the interval $\{x\mid (x,y)\in \Omega\}$. The notation $v$ means that $y$ is the vertical
 component.
\item $I_{\Omega,h}(x):=\{y:(x,y)\in \Omega\}$.
\item $I_{\Omega,v}(y):=\{x:(x,y)\in \Omega\}$.
\item $a\wedge b=\min \{a,b\}$ for $a,b\in \bR$.
\item $a\vee b=\max\{a,b\}$ for $a,b\in \bR$.
\end{itemize}
For clarity, we assume that $\Omega$ is a horizontally and vertically convex domain such that its horizontal cross sections and vertical cross sections are closed intervals. Then $\Omega$ can be written as
\[
\Omega=\{(x,y)\in \bR^2: a\leq x\leq b, g(x)\leq y\leq h(x)\}
\]
and
\[
\Omega=\{(x,y)\in \bR^2: c\leq y\leq d, u(y)\leq x\leq v(y)\},
\]
where $g,h$ and $u,v$ are functions on the intervals $[a,b]$ and $[c,d]$, respectively. Before our discussion of the almost orthogonality estimate,
we shall define some expanded domains for $\Omega$.
\begin{defn}\label{Def expanded domains I}
Assume $\Omega$ is a domain given as above. Then $\Omega_h^{*}$ is said to be a horizontally expanded domain for $\Omega$ if there exists a positive number $\epsilon$ and a nonnegative function $\alpha$ on $[c,d]$ such that
\[
\Omega_h^{*}=\{(x,y):c\leq y\leq d, u(y)-\alpha(y)\leq x\leq v(y)+\alpha(y)\}
\]
and
\[
\delta_{\Omega_h^{*},v}(y)\geq (1+2\epsilon)\delta_{\Omega,v}(y), \quad c\leq y\leq d,
\]
i.e.
\[
v(y)-u(y)+2\alpha(y)\geq (1+2\epsilon)(v(y)-u(y)).
\]
Similarly, we call $\Omega_v^{*}$ a vertical expanded domain for $\Omega$ if there exists an $\epsilon>0$
and a nonnegative function $\beta$ such that
\[
\Omega_v^{*}=\{(x,y):a\leq x\leq b, g(x)-\beta(x)\leq y\leq h(x)+\beta(x)\}
\]
and
\[
\delta_{\Omega_v^{*},h}(x)\geq (1+2\epsilon)\delta_{\Omega,h}(x), \quad x\in[a,b].
\]
\end{defn}
\begin{defn}\label{Def expanded domains II}
Assume $\Omega$ is a horizontally and vertically convex domain. A domain $\Omega^{*}$ is called an expanded domain
for $\Omega$ if there exist two expanded domains $\Omega_{h}^{*}$ and $\Omega_{v}^{*}$ as in Definition {\rm \ref{Def expanded domains I} } such that
\[
\Omega^{*}=\Omega_{h}^{*}\cup\Omega_{v}^{*}.
\]
\end{defn}
Now we shall state the almost orthogonality estimate between two oscillatory integral operators which are supported on horizontally
and vertically convex domains. It should be pointed out that the same estimate had been established by Phong and Stein for oscillatory
integral operators supported on curved trapezoids; see Phong-Stein \cite{PS1998}.
\begin{lemma}\label{lemma almost orthogonal esti}
Assume $T_1$ and $T_2$ are two oscillatory integral operators defined as {\rm(1.1)} with a common phase $S$, where the cut-off functions for $T_1$ and $T_2$ are $\varphi_1$ and $\varphi_2$, respectively.
Suppose $\Omega_1$ and $\Omega_2$ are two horizontally and vertically convex domains such that
\[
\supp(\varphi_1)\subseteq  \Omega_{1},\quad \supp(\varphi_2)\subseteq  \Omega_2.
\]
If all of the following conditions hold,
\begin{enumerate}
\item[{\rm (i)}] For two positive constants $\mu, A>0$, there exist two expanded domains $\Omega_1^{*}$ and
$\Omega_2^{*}$ such that
$\mu\leq \l|S_{xy}''\r|\leq A\mu$ on $\Omega_1^{*}.$
The Hessian $S_{xy}''$ is a polynomial type function with order $ord(S)$ in $y$, uniformly in $x$, on vertical cross sections of appropriately expanded domains $\Omega_1^{\ast}$ and $\Omega_2^{\ast}$.
\item[{\rm (ii)}] For any horizontal line segment $L$ joining a point $(x,z)\in \Omega_1^{*}$ and another one
$(y,z)\in\Omega_2^{*}$, the Hessian $S_{xy}''$ does not change sign on $L$ and $\l|S_{xy}''(x,y)\r|\leq A\mu$ for
all points $(x,y)\in L$.
\item[{\rm (iii)}] There exists a constant $B\geq 1$ such that $I_{\Omega_{2}^{*},v}(y)\subseteq I^{\ast}_{\Omega_{1}^{*},v}(y; B)$
for every $y\in \bR$. Here the notation $I^{\ast}(B)$ denotes the interval concentric with $I$ and having length $|I^{\ast}(B)|=B|I|$.
\item[{\rm (iv)}] There are two constants $M_1,M_2>0$ such that for all $(x,y)\in\Omega_i$
\[
\sum_{k=0}^2\left(\delta_{\Omega_i,h}(x)\right)^k\left|\partial_y^k\varphi_i(x,y)\right|\leq M_i, \quad i=1,2,
\]
\end{enumerate}
\noindent then there exists a constant $C=C(ord(S),C_{S''_{xy}},\epsilon_1,\epsilon_2,A,B)$ such that
\[
\l\|T_1T_2^{*}\r\|\leq CM_1M_2|\lambda\mu|^{-1},
\]
where $\epsilon_1,\epsilon_2>0$ are expanded factors (see Definition {\rm \ref{Def expanded domains I}}) for $\Omega_1^{\ast}$ and $\Omega_2^{\ast}$.
\end{lemma}
\begin{proof}
It is easy to see that the integral kernel associated with $T_1T_2^{*}$ is given by
\[
K(x,y)=\int_{-\infty}^{+\infty}
e^{i\lambda[S(x,z)-S(y,z)]}\varphi_1(x,z)\overline{\varphi_2(y,z)}dz.
\]
Define a function $\Phi$ by
\[
\Phi(z)=\partial_zS(x,z)-\partial_zS(y,z)
\]
and a linear differential operator $\mathcal{D}$ by
\[
\mathcal{D}f(z)=\frac{1}{i\lambda}\frac{1}{\Phi(z)}f'(z).
\]
For convenience, we also define the transpose $\mathcal{D}^t$ of $\mathcal{D}$ as
\[
\mathcal{D}^tf(z)=-\frac{1}{i\lambda}\frac{\partial}{\partial z}\l[\frac{f(z)}{\Phi(z)}\r].
\]
Since $\supp(\varphi_i)\subseteq  \Omega_i$ for $i=1,2$, we see that if $\varphi_1(x,z)\overline{\varphi_2(y,z)}\neq 0$
then $z\in I_{\Omega_1,h}(x)\cap I_{\Omega_2,h}(y)$. By the assumptions (i), (ii) and (iii), there exists a constant
$C=C(B)$ such that
\[
|\Phi(z)|=\left|\int_{y}^x\partial_z\partial_uS(u,z)du\right|
\geq
C\mu|x-y|, \quad z\in I_{\Omega_1,h}(x)\cap I_{\Omega_2,h}(y).\]
Combining this result, we deduce from the assumption (ii) that
\[
C\mu|x-y|\leq |\Phi(z)|\leq A\mu|x-y|, \quad z\in I_{\Omega_1^{*},h}(x)\cap I_{\Omega_2^{*},h}(y).
\]
As assumed above, $I_{\Omega_1,h}(x)\cap I_{\Omega_2,h}(y)\neq \varnothing$. Choose an arbitrary point $a$ from this intersection. By Definitions \ref{Def expanded domains I} and \ref{Def expanded domains II}, we see that
\begin{equation}\label{section 2 nonempty intersection property I}
\Big(a-\epsilon_1\delta_{\Omega_1,h}(x), a+\epsilon_1 \delta_{\Omega_1,h}(x)\Big)
\subseteq I_{\Omega_1^{*},h}(x),~~~
\Big(a-\epsilon_2\delta_{\Omega_2,h}(y), a+\epsilon_2
\delta_{\Omega_2,h}(y)\Big)
\subseteq I_{\Omega_2^{*},h}(y).
\end{equation}
Hence
\begin{eqnarray}\label{section 2 nonempty intersection property II}
|I_{\Omega_1^{*},h}(x)\cap I_{\Omega_2^{*},h}(y)|
&\geq&
2\min\big\{\epsilon_1\delta_{\Omega_1,h}(x),
\epsilon_2 \delta_{\Omega_2,h}(y)\big\}\nonumber\\
&\geq&
2\epsilon_1 \wedge \epsilon_2
|I_{\Omega_1,h}(x)| \wedge |I_{\Omega_2,h}(y)|,
\end{eqnarray}
where $\epsilon_1,\epsilon_2$ appear in the definition of $\Omega_1^{*}$ and $\Omega_2^{*}$. By the scaling property
of polynomials in Lemma \ref{sec2 esti poly tpype func}, there exists a constants $C$ depending only on $ord(S)$ and $C_{S''_{xy}}$ such that
\begin{align*}
\sup\l|\partial_z^k\Phi(z)\r|&\leq C\mu|x-y|\l|I_{\Omega_1^{*},h}(x)\cap I_{\Omega_{2}^{*},h}(y)\r|^{-k}\\
&\leq C(\epsilon_1\wedge\epsilon_2)^{-k}\mu|x-y|
\left( \l|I_{\Omega_1^{*},h}(x)\r|\wedge\l|I_{\Omega_2^{*},h}(y)\r|
\right)^{-k}.
\end{align*}
By integration by parts, we have
\begin{align*}
K(x,y)&=\int_{-\infty}^{+\infty}
\mathcal{D}^2\left(e^{i\lambda[S(x,z)-S(y,z)]}\right)
\varphi_1(x,z)\overline{\varphi_2(y,z)}dz\\
&=\int_{-\infty}^{+\infty}
e^{i\lambda[S(x,z)-S(y,z)]}
\l(\mathcal{D}^t\r)^2\left(\varphi_1(x,z)\overline{\varphi_2(y,z)}\right)dz.
\end{align*}
It is easy to see that $\l(\mathcal{D}^t\r)^2\left(\varphi_1(x,z)\overline{\varphi_2(y,z)}\right)$ is a linear combination of the following terms
\[
\frac{1}{(i\lambda)^2}\frac{\partial^{k_1}}{\partial z^{k_1}}\left(\frac{1}{\Phi(z)}\right)\frac{\partial^{k_2}}{\partial z^{k_2}}
\left(\frac{1}{\Phi(z)}\right)\frac{\partial^{k_3}}{\partial z^{k_3}}\varphi_1(x,z)\frac{\partial^{k_4}}{\partial z^{k_4}}\overline{\varphi(y,z)}
\]
with $k_i$ being nonnegative integers satisfying
$k_1+k_2+k_3+k_4=2.$
By induction, we can prove that
\begin{align*}
\left|\frac{d^k}{dz^k}\left(\frac{1}{\Phi(z)}\right)\right|&\leq C(ord(S),C_{S''_{xy}}, A)\left(\mu|x-y|\right)^{-1}
\l|I_{\Omega_1^{*},h}(x)\cap I_{\Omega_2^{*},h}(y)\r|^{-k}\\
&\leq C(ord(S),C_{S''_{xy}},A,\epsilon_1,\epsilon_2)(\mu|x-y|)^{-1}
\big( \l|I_{\Omega_1,h}(x)\r|\wedge \l|I_{\Omega_2,h}(y)\r| \big)^{-k}.
\end{align*}
In fact, by Remark \ref{Rmk polynomial type func}, we have claimed that $(\Phi(z))^{-1}$ satisfies the estimate in Lemma \ref{sec2 esti poly tpype func} on the interval $I_{\Omega_1^{*},h}(x)\cap I_{\Omega_2^{*},h}(y)$. By the assumption (iv), it is clear that
\[
\left|\partial_z^k\varphi_i(u,z)\right|\leq M_i\left(\l|I_{\Omega_1,h}(x)\r|\wedge \l|I_{\Omega_2,h}(y)\r|\right)^{-k}
\]
for $u\in\{x,y\}$ and $z\in I_{\Omega_1,h}(x)\cap I_{\Omega_2,h}(y)$. Combining all above estimates, we obtain
\[
\left|\left(\mathcal{D}^t\right)\left(\varphi_1(x,z)\overline{\varphi_2(y,z)}\right)\right|\leq CM_1M_2
(\mu\lambda|x-y|)^{-2}\left(\l|I_{\Omega_1,h}(x)\r|\wedge \l|I_{\Omega_2,h}(y)\r|\right)^{-2}.
\]
On the other hand, if we take absolute value into the integral for $K$, then we find that
\[
|K(x,y)|\leq CM_1M_2\l|I_{\Omega_1,h}(x)\r|\wedge \l|I_{\Omega_2,h}(y)\r|.
\]
Hence we obtain
\[
|K(x,y)|\leq CM_1M_2\left(1+\delta_{\Omega_1,h}(x)\wedge \delta_{\Omega_2,h}(y)|\lambda|\mu|x-y|\right)^{-2}\cdot
\delta_{\Omega_1,h}(x)\wedge \delta_{\Omega_2,h}(y).
\]
Recall that $\delta_{\Omega_1,h}(x)=\l|I_{\Omega_1,h}(x)\r|$ and $\delta_{\Omega_2,h}(y)=\l|I_{\Omega_2,h}(y)\r|$.
Set $a=\delta_{\Omega_1,h}(x)$ and $b=\delta_{\Omega_2,h}(y)$. We have
\begin{align*}
&\int_{\bR^2}|K(x,y)||f(y)||g(x)|dxdy\\
&\leq CM_1M_2
\left(\int_{\bR^2}\frac{a}{(1+a\mu|\lambda||x-y|)^2 }|f(y)||g(x)|dxdy+\int_{\bR^2}\frac{b}{(1+b\mu|\lambda||x-y|)^2}|f(y)||g(x)|dxdy \right)\\
&\leq CM_1M_2|\lambda\mu|^{-1}\int_{\bR}|g(x)|\mathcal{M}f(x)dx+CM_1M_2|\lambda\mu|^{-1}\int_{\bR}|f(y)|\mathcal{M}g(y)dy\\
&\leq CM_1M_2|\lambda\mu|^{-1}\|f\|_{L^2}\|g\|_{L^2},
\end{align*}
where $\mathcal{M}$ is the Hardy-Littlewood maximal operator and $C=C(ord(S),C_{S''_{xy}},\epsilon_1,\epsilon_2,A,B)$. Here we have used the fact that $\mathcal{M}$ is a bounded operator on $L^2$. The proof of the lemma is complete.
\end{proof}

Now we shall give an interpolation lemma with change of power weights. One can compare this lemma with the well-known Stein-Weiss interpolation with change of measures; see Stein-Weiss \cite{stein-weiss}. The following lemma is an extension of Proposition 1.4 in \cite{pansampson}; see \cite{ShiYan} and \cite{Shi} for its proof.
\begin{lemma}\label{section 2 stein-weiss interpolation}
Let $dx$ be the Lebesgue measure on $\bR$. Assume $T$ is a linear operator defined on all simple functions with respect to $dx$. If there exist
two constants $A_1,A_2>0$ such that
\begin{enumerate}
\item[{\rm (i)}]  $\|Tf\|_{L^{\infty}(dx)}\leq A_1\|f\|_{L^1(dx)}$ for all simple functions $f$;
\item[{\rm (ii)}] $\||x|^aTf\|_{L^{p_0}(dx)}\leq A_2\|f\|_{L^{p_0}(dx)}$ for some $1<p_0<+\infty$, $a\in \bR$
satisfying $ap_0\neq -1$;
\end{enumerate}
then for any $\theta\in (0,1)$ there exists a constant $C=C(a,p_0,\theta)$ such that
\[
\||x|^bTf\|_{L^p(dx)}\leq CA_1^{\theta}A_2^{1-\theta}\|f\|_{L^p(dx)}
\]
for all simple functions $f$, where $b$ and $p$ satisfy
$
b=-\theta+(1-\theta)a$ and $\frac{1}{p}=\theta+\frac{1-\theta}{p_0}.
$

\end{lemma}

We shall need a uniform estimate with non-sharp decay exponent in Section 4.
\begin{lemma}\label{Lemma uniform L2 estimate}
Assume $S$ is a real-valued smooth function in the unit square $Q=(0,1)^2$. Let $T_\lambda$ be the oscillatory integral operator in {\rm (\ref{section 1 def Tlambda})} with phase
$S$ and cut-off function $\varphi\in C^{\infty}_0(Q)$. If there exist two positive integers $j$ and $k$ such that $\partial_x^j\partial_y^k S(x,y)\geq 1$ on $Q$, then we have the following results.\\

\noindent {\rm (i)} In the case $j=k=1$, if in addition there exists a positive integer $N$ such that either $\partial_x^{N+1}\partial_y S$ or $\partial_x\partial_y^{N+1} S$ is single-signed in $Q$, then we have
\begin{equation}\label{Uniform estimates for j=k=1}
\|T_\lambda f\|_{L^2}\leq C(N,\varphi)|\lambda|^{-1/2}
\log^{1/2}(2+|\lambda|)\|f\|_{L^2}.
\end{equation}

\noindent {\rm (ii)} If either $j\geq 2$ or $k\geq 2$, then there exists a constant $C=C(j,k,\varphi)$ such that
\begin{equation}\label{Uniform estimates for j,k greater than 1}
\|T_\lambda f\|_{L^2}\leq C|\lambda|^{-\delta}
\|f\|_{L^2},~~~~~\delta=\frac{1}{2}\cdot\frac{1}{(j\vee k)2^{j\wedge k}}.
\end{equation}
\end{lemma}

For the proof of this lemma, we refer the reader to Carbery-Christ-Wright \cite{carbery} where the optimal decay exponent $\delta=\frac{1}{2(j\vee k)}$ was obtained when either $j$ or $k$ is equal to $1$ and at least one of them is greater than $1$. In the same paper \cite{carbery}, uniform estimates were also given for general $j\geq 2$ or $k\geq 2$ with non-sharp decay exponents. Especially, if $S$ is a polynomial, one can see Ricci-Stein \cite{ricci} for $\delta<\frac{1}{2deg(S)}$ and Phong-Stein-Sturm \cite{PSS2001} for the general optimal decay exponent $\delta=\frac{1}{2(j\vee k)}$. For further related results under more general nondegenerate conditions, one can see Christ-Li-Tao-Thiele \cite{Christ-Li-Tao-Thiele};
see also Greenblatt \cite{greenblatt3} for uniform multilinear decay estimates in the case of smooth phases.

Now we state the operator van der Corput lemma for nondegenerate oscillatory integral operators supported on horizontally and vertically convex domains. This lemma was established by Phong-Stein-Sturm \cite{PSS2001}. For its proof, one can also see \cite{PSS2001}.
\begin{lemma}\label{Lemma operator vers van der Corput}
Let $T_\lambda$ be defined as {\rm (\ref{section 1 def Tlambda})} with the cut-off $\varphi$ supported in a curved trapezoid $\Omega$. Assume the Hessian $S_{xy}''$ is a polynomial type function with order $ord(S)$ in $y$, uniformly in $x$, on each vertical cross section $I_{\Omega,h}(x)=\{y:(x,y)\in\Omega\}$. If the following assumptions are ture:
\begin{enumerate}
\item[{\rm (i)}] For some $\mu, A>0$, the Hessian $S_{xy}''$ satisfies
$\mu\leq\l|S_{xy}''(x,y)\r|\leq A\mu$ for $(x,y)\in \Omega;$
\item[{\rm (ii)}] Let $\delta_{\Omega,h}(x)$ be the length of $I_{\Omega,h}(x)$. For some $M>0$,
\[
\sum_{k=0}^2\sup_{\Omega}\left(\delta_{\Omega,h}(x)\right)^k
|\partial_y^k\varphi(x,y)|\leq M,
\]
\end{enumerate}
\noindent then there exists a constant $C=C(ord(S),C_{S''_{xy}},A)$ such that for all $f\in L^2(\bR)$
\[
\|T_\lambda f\|_{L^2}\leq CM(|\lambda|\mu)^{-1/2}\|f\|_{L^2}.
\]
\end{lemma}

In Section 3, the following lemma will be frequently used to balance oscillation and size estimates of nondegenerate oscillatory integral operators.
\begin{lemma}\label{section 2 control of abs sum}
For arbitrary $0<\epsilon<1$ and $M>1$, there exists a constant $C=C(\epsilon, M)>0$ such that for all $A, B>0$,
\begin{equation*}
\sum_{k\in\bZ}\min\left\{A\epsilon^k, BM^k\right\}\leq C\sup_{k\in \bZ}\left(\min\left\{A\epsilon^k, BM^k\right\}\right).
\end{equation*}
As a consequence, if $\theta\in(0, 1)$ satisfies $\epsilon^{\theta}M^{1-\theta}=1$, then
\[
\sum_{k\in\bZ}\min\left\{A\epsilon^k, BM^k\right\}\leq CA^{\theta}B^{1-\theta}.
\]
\end{lemma}
\begin{proof}
Since $\{A\epsilon^k\}$ is strictly decreasing and $\{BM^k\}$ is strictly increasing, there exists a unique $k_0\in\bZ$ such that
$A\epsilon^{k_0}\geq BM^{k_0}$ and $A\epsilon^{k_0+1}<BM^{k_0+1}.$
For this $k_0$, it is easy to see that $A\epsilon^{k_0}$ and $BM^{k_0}$ are comparable up to constants depending only on $\epsilon$ and $M$. Now we have
\begin{align*}
\sum_{k\in\bZ}\min\l\{A\epsilon^k, BM^k\r\}&=\sum_{k\leq k_0}BM^k+\sum_{k>k_0}A\epsilon^k\\
&\leq C(\epsilon,M)\min\left\{A\epsilon^{k_0}, BM^{k_0}\right\}\\
&\leq CA^{\theta}B^{1-\theta}
\end{align*}
provided that $\epsilon^{\theta}M^{1-\theta}=1$. The proof is complete.
\end{proof}

Now we give the well-known Puiseux factorization theorem for real-analytic functions (also for formal power series in $\bR^2$). Assume $n$ is a positive integer. Let $\mathbb{C}[[{x^{1/n}}]]$ be the ring of formal fractional power series in $x^{1/n}$ with complex coefficients. For each $r\in \mathbb{C}[[{x^{1/n}}]]$, i.e., $r(x)=\sum_{i\geq 0}a_ix^{i/n}$, the series $\sigma_{\epsilon}(r)(x)=\sum_{i\geq 0}\epsilon^ia_ix^{i/n}$ will be called the conjugate series of $s$ for all $n$-th roots $\epsilon$ of the unit. For the concept of Puiseux series and its polydromy order $n$, one can see Casas-Alvero {\rm\cite{casas-alvero}}.
\begin{theorem}\label{Puiseux theorem}
{\rm(\cite{casas-alvero}, p32, Theorem 1.8.3)} Let $S$ be a real-analytic function near the origin in $\bR^2$. Then there exists a non-vanishing real-analytic function $U$ and finitely many Puiseux series $r_{\nu}(x)$, both defined in a small neighborhood of the origin, such that
\begin{equation}
S(x,y)=U(x,y)\prod_{\nu}(y-r_{\nu}(x)).
\end{equation}
Moreover, for each root $r_{\nu}$, all of its conjugate series also appear in the above decomposition with the same multiplicity.
\end{theorem}

\section{$L^2$ Damping Decay Estimates}
Assume $S$ is a real-valued analytic function near the origin in $\mathbb{R}^2$. Now we only consider the degenerate case $S_{xy}''(0,0)= 0$. By Taylor's expansion, we can write $S(x,y)=\sum_{k,l\geq 0} a_{k,l}x^ky^l$ in a small neighborhood of the origin. Then the Newton polyhedron $N(S)$ is defined by
\begin{equation}
N(S):={\rm Convex~hull~of~}\Big( \bigcup_{a_{k,l}\neq 0}\Big\{(x,y):x\geq k, y\geq l\Big\} \Big)
\end{equation}
with the union taken over all $k,l\geq 0$ for which $a_{k,l}\neq 0$.

By the Puiseux expansion theorem (see Casas-Alvero \cite{casas-alvero} for example), the Hessian $S_{xy}''$ can be decomposed as the following form:
\begin{equation}\label{Section 3 expansion of the Hessin}
S_{xy}''(x,y)=U(x,y)x^my^s\prod_{\nu}(y-r_\nu(x)),
\end{equation}
where $U$ is a nonvanishing smooth factor near the origin, $m$ and $s$ are nonnegative integers, and $r_\nu$ are nontrivial roots which can be written as Puiseux series. To establish the desired damping estimates, we need a complete classification of the roots $r_{\nu}$. The classification method due to Phong-Stein \cite{PS1997} will be used in this section.

We order the leading exponents of all roots $r_\nu$ as follows:
\[0<a_1<\cdots<a_l<a_{l+1}<\cdots<a_n<+\infty.\]
Denote by
$\Psi\begin{bmatrix}
\begin{smallmatrix}
\cdot\\l
\end{smallmatrix}
\end{bmatrix}
$
the cluster of nonzero roots $r_{\nu}$ in (\ref{Section 3 expansion of the Hessin}) of form
\begin{equation}\label{Section 3 Expansion for nontrivila roots}
r_\nu(x)=C_{l}^{\alpha}x^{a_l}+o(x^{a_l}),~~~x\rightarrow 0,
\end{equation}
where $C_l^{\alpha}$ is a nonzero complex number. In this cluster,
each root is repeated according to its multiplicity. Hence the number of nontrivial roots is counted with the multiplicity of each root.

To distinguish the roots with the same leading exponent $a_l$, we consider the set of all distinct coefficients $C_l^\alpha$ of $x^{a_l}$. Here $l$ is the index for the leading exponents $a_l$, and $\alpha$ is the index for the coefficients of $x^{a_l}$. Then we can see
\[\Psi
\begin{bmatrix}
\cdot\\
l
\end{bmatrix}=\bigcup_{\alpha}\Psi
\begin{bmatrix}
\alpha\\
l
\end{bmatrix},\]
where each root in
$
\Psi
\begin{bmatrix}
\begin{smallmatrix}
\alpha\\
l
\end{smallmatrix}
\end{bmatrix}$ is of the form (\ref{Section 3 Expansion for nontrivila roots}) for fixed $C_l^{\alpha}$. If we continue this process, a complete classification of all roots can be obtained. Indeed, we have
\[\Psi
\begin{bmatrix}
\alpha_1,&\cdots,&\alpha_{N-1},&\cdot\\
l_1,&\cdots,&l_{N-1},&l_N
\end{bmatrix}=\bigcup_{\alpha_N}\Psi
\begin{bmatrix}
\alpha_1,&\cdots,&\alpha_{N-1},&\alpha_N\\
l_1,&\cdots,&l_{N-1},&l_N
\end{bmatrix},\]
where the cluster $\Psi\begin{bmatrix}
\alpha_1,&\cdots,&\alpha_{N-1},&\alpha_N\\
l_1,&\cdots,&l_{N-1},&l_N
\end{bmatrix}$
consists of all roots $r_\nu$ of form
\[r_\nu(x)=C_{l_1}^{\alpha_1} x^{a_{l_1}}+C_{l_1,l_2}^{\alpha_1,\alpha_2}x^{a_{l_1,l_2}^{\alpha_1}}+\cdots+
C_{l_1,\cdots,l_N}^{\alpha_1,\cdots,\alpha_N}x^{a_{l_1,\cdots,l_N}^{\alpha_1,\cdots,\alpha_{N-1}}}+
o(x^{a_{l_1,\cdots,l_N}^{\alpha_1,\cdots,\alpha_{N-1}}}).\]
To make precisely the number of nonzero roots in the above clusters, we use
the following notations:
\[\textbf{N}\begin{bmatrix}
\alpha_1,&\cdots,&\alpha_{k-1},&\cdot\\
l_1,&\cdots,&l_{k-1},&l_k
\end{bmatrix}=\#\ \text{roots in}\ \Psi\begin{bmatrix}
\alpha_1,&\cdots,&\alpha_{k-1},&\cdot\\
l_1,&\cdots,&l_{k-1},&l_k
\end{bmatrix},\]
\[\textbf{N}\begin{bmatrix}
\alpha_1,&\cdots,&\alpha_{k-1},&\alpha_k\\
l_1,&\cdots,&l_{k-1},&l_k
\end{bmatrix}=\#\ \text{roots in}\ \Psi\begin{bmatrix}
\alpha_1,&\cdots,&\alpha_{k-1},&\alpha_k\\
l_1,&\cdots,&l_{k-1},&l_k
\end{bmatrix}.\]
Here $\#A$ means the cardinality of the set $A$. Sometimes, we also use the notation $N$ to denote the number of nontrivial roots of $S_{xy}''$, counting with multiplicities. In other words,
\begin{equation}\label{Section 3 total number of roots}
N=\sum_{l=1}^{n}\textbf{N}\begin{bmatrix}
\cdot\\
l
\end{bmatrix}=\sum_{l=1}^\textbf{N}\#\ \text{roots in}\ \Psi\begin{bmatrix}
\cdot\\
l
\end{bmatrix}.
\end{equation}
Assume $(A_0, B_0), (A_1, B_1), \cdots, (A_{n}, B_{n})$ are vertices of the Newton polyhedron of $S_{xy}''$. For clarity, we also assume $A_0<A_1<\cdots<A_{n}$ and $B_0>B_1>\cdots>B_{n}$.
One can verify that
\begin{equation}
A_0=m,~~~~~ B_0=s+\sum_{l=1}^n \textbf{N}\begin{bmatrix}
\cdot\\
l
\end{bmatrix},
\end{equation}
and for each $1\leq r \leq n$,
\begin{equation}
A_r=m+\sum_{l=1}^r\textbf{N}\begin{bmatrix}
\cdot\\
l
\end{bmatrix}a_l,~~~~~
 B_r=s+\sum_{l=r+1}^n \textbf{N}\begin{bmatrix}
\cdot\\
l
\end{bmatrix}.
\end{equation}

With the above notations, we can give an equivalent formulation of Theorem \ref{section 1 Main theorem} in terms of the vertices of $N(S)$.
\begin{theorem}\label{section 3 Thm main result}
Assume $S$ is a real-valued analytic function near the origin and $T_\lambda$ is defined as {\rm (\ref{section 1 def Tlambda})}. Then there exists a constant $C$, independent of $\lambda$ and $f$, such that
\begin{equation}\label{sec-main-est}
\|T_\lambda f\|_{L^{p_r}}\leq C|\lambda|^{-\frac{1}{A_r+B_r+2}}\|f\|_{L^{p_r}}, \quad p_r=\frac{A_r+B_r+2}{A_r+1}
\end{equation}
for all $0\leq r\leq n$.
\end{theorem}
To prove the estimate \eqref{sec-main-est}, we shall consider a family of damped oscillatory integral operators
\begin{equation}\label{damp-opr}
W_zf(x)=\int_{-\infty}^{+\infty}e^{i\lambda S(x,y)}|D(x,y)|^z\varphi(x,y)f(y)dy
\end{equation}
with the damping factor $D$ defined by
\begin{equation}\label{section 3 Def damping D}
D(x,y)=x^m\prod_{l=1}^r\prod_{r_\nu\in \Psi\begin{bmatrix}
\begin{smallmatrix}
\cdot\\
l
\end{smallmatrix} \end{bmatrix}}\left(y-r_\nu(x)\right).
\end{equation}
Then we have the following $L^2$ decay estimate.
\begin{theorem}\label{Thm L2 damping decay estiamte}
Assume $W_z$ is defined as above. For $z\in\mathbb{C}$ with real part
\begin{equation}
\Re(z) = \frac{A_r-B_r}{2A_r}\cdot\frac{1}{1+B_r},
\end{equation}
there exists a constant $C=C(S,\varphi)$ such that
\begin{equation}\label{Section 3 main L2 decay estimate}
\|W_zf\|_{L^2}~
\leq C(1+|z|^2)|\lambda|^{-\sigma}\|f\|_{L^2}~~~\textrm{with}~~~
\sigma=\frac{1}{2(1+B_r)}.
\end{equation}
\end{theorem}
\begin{remark}\label{sec3 Remark}
By Theorem {\rm\ref{Puiseux theorem}}, each product $\prod_{r_\nu\in \Psi\begin{bmatrix}
\begin{smallmatrix}
\cdot\\
l
\end{smallmatrix} \end{bmatrix}}\left(y-r_\nu(x)\right)$ in {\rm(\ref{section 3 Def damping D})} is conjugation invariant and hence $D$ is a real-analytic function
{\rm(Lemma 1.2.1 in {\rm\cite{casas-alvero}})}. Since $S$ is real valued, $S(x,r_{\nu}(x))\equiv 0$ if and only if $S(x,\overline{r_{\nu}}(x))\equiv 0$. This implies that $D$ is also real valued.
\end{remark}
\begin{proof}
Choose a cut-off function $\Phi\in C_{0}^\infty$ such that
$\supp(\Phi)\subseteq [\frac{1}{2}, 2]$ and $\sum_{j\in \mathbb{Z}}\Phi\left(\frac{x}{2^j}\right)=1$ for all $x>0$. Write $W_z=\sum_{\sigma_1,\sigma_2=\pm}W_z^{\sigma_1,\sigma_2}$
and
$W_z^{\sigma_1,\sigma_2}
=\sum_{j,k}W_{j,k}^{\sigma_1,\sigma_2}$
with
\begin{equation}\label{damp-decom}
W_{j,k}^{\sigma_1,\sigma_2}f(x)
=
\int_{-\infty}^{+\infty}e^{i\lambda S(x,y)}|D(x,y)|^z\Phi\left(\sigma_1\frac{x}{2^j}\right)
\Phi\left(\sigma_2\frac{y}{2^k}\right)\varphi(x,y)f(y)dy.
\end{equation}
In what follows, it will be shown that the estimate (\ref{Section 3 main L2 decay estimate}) is true with $W_z$ replaced by $W_z^{+,+}$. Other operators $W_z^{\sigma_1,\sigma_2}$ can be treated similarly. For convenience, we shall still write $W_z$ instead of $W_z^{++}$.\\

\textbf{Case (i)}\quad $k\gg a_1j.$

\noindent The notation $k\gg a_1j$ means that $k\geq a_1j+N_0$ for some large $N_0>0.$ Then both the Hessian $S_{xy}''$ and $D$ have sizes equivalent to two positive constants, up to a multiplicative constant.

Assume first $m>0$. Since $W_{j,k}W_{j',k'}^{*}=0$ for $|k-k'|\geq 2$, we assume $|k-k'|\leq 1$. Then by Lemma \ref{lemma almost orthogonal esti},
\begin{equation*}
\left\|W_{j,k}W_{j',k'}^{*}\right\|\leq C\left[|\lambda|2^{mj}2^{k(N+s)}\right]^{-1}
\Big[2^{mj}2^{k\cdot\sum\limits_{l=1}^r\textbf{N}\begin{bmatrix}
\begin{smallmatrix}
\cdot\\
l
\end{smallmatrix}
\end{bmatrix}}\Big]^{\Re(z)}
\Big[2^{mj'}2^{k'\cdot\sum\limits_{l=1}^r\textbf{N}\begin{bmatrix}
\begin{smallmatrix}
\cdot\\
l
\end{smallmatrix}
\end{bmatrix}}\Big]^{\Re(z)}
~~~j\geq j',
\end{equation*}
where $N$ is given by (\ref{Section 3 total number of roots}).

By the size estimate, the $L^2$ operator norm $\|W_{j,k}^{}W_{j',k'}^{*}\|$
is bounded by
\[ \left\|W_{j,k}\right\|\left\|W_{j',k'}^{*}\right\|\leq C\Big[2^{mj}2^{k\cdot\sum\limits_{l=1}^r\textbf{N}\begin{bmatrix}
\begin{smallmatrix}
\cdot\\
l
\end{smallmatrix}
\end{bmatrix}}\Big]^{\Re(z)}
\Big[2^{mj'}2^{k'\cdot\sum\limits_{l=1}^r\textbf{N}\begin{bmatrix}
\begin{smallmatrix}
\cdot\\
l
\end{smallmatrix}
\end{bmatrix}}\Big]^{\Re(z)} 2^{j/2}2^{k/2}2^{j'/2}2^{k'/2}.
\]
Put $\theta=(1+B_r)^{-1}.$ Taking a convex combination of the above two estimates, we see that $\left\|W_{j,k}W_{j',k'}^{*}\right\|$ is bounded by a constant multiple of
\begin{eqnarray*}
& &\left[|\lambda|2^{mj}2^{k(N+s)}\right]^{\theta}
\Big[2^{mj}2^{k\cdot\sum\limits_{l=1}^r\textbf{N}\begin{bmatrix}
\begin{smallmatrix}
\cdot\\
l
\end{smallmatrix}
\end{bmatrix}}\Big]^{\Re(z)}
\Big[2^{mj'}2^{k'\cdot\sum\limits_{l=1}^r\textbf{N}\begin{bmatrix}
\begin{smallmatrix}
\cdot\\
l
\end{smallmatrix}
\end{bmatrix}}\Big]^{\Re(z)}\Big(2^{j/2}2^{k/2}\Big)^{1-\theta}
\Big(2^{j'/2}2^{k'/2}\Big)^{1-\theta}\\
&\leq & C|\lambda|^{-\theta}2^{-|j-j'|\gamma}, \quad \gamma=m\Re(z)+\frac{1-\theta}{2}.
\end{eqnarray*}
By direct calculation, we see that the sum of the exponents of $2^k$ and $2^{k'}$ is equal to, using $2\Re(z)=\left(1-\frac{B_r}{A_r}\right)\theta$ and $1-\theta=B_r\theta$,
\begin{eqnarray*}
-(N+s)\theta+2\sum_{l=1}^r\textbf{N}\begin{bmatrix}
\cdot\\
l
\end{bmatrix}\Re(z)+1-\theta
&=&-(N+s)\theta+\sum_{l=1}^r\textbf{N}\begin{bmatrix}
\cdot\\
l
\end{bmatrix}\left(1-\frac{B_r}{A_r}\right)\theta+B_r\theta\\
&=&-\sum_{l=1}^r\textbf{N}\begin{bmatrix}
\cdot\\
l
\end{bmatrix}\frac{B_r}{A_r}\theta\leq 0.
\end{eqnarray*}
In view of $k\gg a_1j$, we have
\begin{align*}
2^{-k\sum\limits_{l=1}^r\textbf{N}\begin{bmatrix}
\begin{smallmatrix}
\cdot\\
l
\end{smallmatrix}
\end{bmatrix}B_r\theta/A_r}&\leq C2^{-a_1j\sum\limits_{l=1}^r\textbf{N}\begin{bmatrix}
\begin{smallmatrix}
\cdot\\
l
\end{smallmatrix}
\end{bmatrix}B_r\theta/A_r}
\leq C2^{-j(A_r-m)B_r\theta/A_r}.
\end{align*}
If we add this new exponent and the original one of $2^j$, we obtain the resulting exponent
\begin{align*}
&-m\theta+m\Re(z)+\frac{1-\theta}{2}-(A_r-m)B_r\theta/A_r\\
&=-m\theta+m(1-B_r/A_r)\theta/2+B_r\theta/2-(B_r-B_rm/A_r)\theta\\
&=-m\Re(z)-(1-\theta)/2.
\end{align*}
Thus we obtain the desired almost orthogonality estimate. It should be pointed out that the exponent $\gamma:=m\Re(z)+(1-\theta)/2$ is positive. Indeed if $\gamma=0$ then $A_r=B_r=0$ which implies that $T_{\lambda}$ is nondegenerate. A similar estimate is also true for $W_{j,k}^{\ast}W_{j',k'}$.

Now we consider the case $m=0$. We shall study $W_k=\sum_{j}W_{j,k}$ with the summation taken over all $j$ satisfying $a_1j\ll k$. The treatment is somewhat different depending on whether $D(x,y)$ is equal to $S_{xy}''$ up to a non-vanishing factor $U$ in (\ref{Section 3 expansion of the Hessin}). Assume first $S_{xy}''$ is not equal to the Hessian $U(x,y)D(x,y)$. Then by the operator van der Corput Lemma \ref{Lemma operator vers van der Corput}, we have
\[
\|W_{k}\|\leq C\left[|\lambda|2^{k\cdot\sum\limits_{l=1}^n\textbf{N}\begin{bmatrix}
\begin{smallmatrix}
\cdot\\
l
\end{smallmatrix}
\end{bmatrix}} 2^{ks}\right]^{-1/2}\left(2^{k\cdot\sum\limits_{l=1}^r\textbf{N}
\begin{bmatrix}
\begin{smallmatrix}
\cdot\\
l
\end{smallmatrix}
\end{bmatrix}}\right)^{\Re(z)}.
\]
The size estimate for $W_k$ is
\[
\|W_k\|\leq C\left(2^{k\cdot\sum\limits_{l=1}^r\textbf{N}\begin{bmatrix}
\begin{smallmatrix}
\cdot\\
l
\end{smallmatrix}
\end{bmatrix}}\right)^{\Re(z)}2^{k/2}2^{k/(2a_1)}.
\]
In the oscillation estimate, the exponent of $2^k$ is
\begin{align*}
&-\left(s+N\right)\big/2+\left(\sum\limits_{l=1}^r\textbf{N}
\begin{bmatrix}
\cdot\\
l
\end{bmatrix}\right)\Re(z)\\
&=-\sum\limits_{l=1}^r\textbf{N}\begin{bmatrix}
\cdot\\
l
\end{bmatrix}\bigg/2-B_r/2+\sum\limits_{l=1}^r\textbf{N}\begin{bmatrix}
\cdot\\
l
\end{bmatrix}\left(\frac{1}{2}-\frac{B_r}{2A_r}\right)\frac{1}{1+B_r}\\
&=-\sum\limits_{l=1}^r\textbf{N}\begin{bmatrix}
\cdot\\
l
\end{bmatrix}\frac{B_r}{2(1+B_r)}
-\frac{B_r}{2}
-\sum\limits_{l=1}^r\textbf{N}\begin{bmatrix}
\cdot\\
l
\end{bmatrix}\frac{1}{2A_r}\cdot\frac{B_r}{1+B_r}<0
\end{align*}
since the fact $S_{xy}''(x,y)\neq U(x,y)D(x,y)$ implies $B_r>0$. It is clear that the exponent of $2^k$ is positive in the above size estimate. Without loss of generality, we may assume $\supp(\varphi)\subset[-1/2, 1/2]^2$. It follows that $W_{j,k}=0$ if $j\geq0$ or $k\geq 0$. For this reason, we may assume $k\leq 0$ from now on. Replacing the exponent $\frac{1}{a_1}$ by the smaller one $\sum\limits_{l=1}^r\textbf{N}\begin{bmatrix}
\begin{smallmatrix}
\cdot\\
l
\end{smallmatrix}
\end{bmatrix}/A_r$, we obtain the following size estimate,
\[
\|W_{k}\|\leq C2^{k\sum\limits_{l=1}^r\textbf{N}\begin{bmatrix}
\begin{smallmatrix}
\cdot\\
l
\end{smallmatrix}
\end{bmatrix}\Re(z)}
2^{k/2}2^{\frac{k}{A_r}\sum\limits_{l=1}^r\textbf{N}\begin{bmatrix}
\begin{smallmatrix}
\cdot\\
l
\end{smallmatrix}
\end{bmatrix}}.
\]
Taking the absolute value into the summation $\sum\limits_{k\leq 0}W_{k}$, we can deduce from Lemma \ref{section 2 control of abs sum} that
\[
\sum_{k\leq 0}\|W_{k}\|\leq C|\lambda|^{-\frac{1}{2(1+B_r)}}.
\]
Here the convex combination, with $\theta=\frac{1}{1+B_r}$, of the above oscillation and size estimates annihilates the exponent of $2^k$. In fact, we see that
\begin{align*}
&-\theta(s+N)/2+\left(\sum\limits_{l=1}^r\textbf{N}\begin{bmatrix}
\cdot\\
l
\end{bmatrix}\right)\Re(z)
+\frac{1}{2}\left(1+\frac{1}{A_r}\sum\limits_{l=1}^r\textbf{N}
\begin{bmatrix}
\cdot\\
l
\end{bmatrix}\right)(1-\theta)\\
&=-\left(\sum\limits_{l=1}^r\textbf{N}\begin{bmatrix}
\cdot\\
l
\end{bmatrix}\right)\cdot\frac{1}{2(1+B_r)}-\frac{B_r}{2(1+B_r)}
+\left(\sum\limits_{l=1}^r\textbf{N}\begin{bmatrix}
\cdot\\
l
\end{bmatrix}\right)
\left(\frac{1}{2}-\frac{B_r}{2A_r}\right)\cdot\frac{1}{1+B_r}+\\ &~~~\frac{1}{2}\left(1+\frac{1}{A_r}\sum\limits_{l=1}^r\textbf{N}
\begin{bmatrix}
\cdot\\
l
\end{bmatrix}\right)\cdot\frac{B_r}{1+B_r}=0.
\end{align*}
It remains to consider the special case $S_{xy}''(x,y)=U(x,y)D(x,y)$ for some nonvanishing factor $U$. The optimal decay $|\lambda|^{-1/2}$ had been established by Phong-Stein \cite{PS1998}. In this case, we must explore the almost orthogonality between $W_k$ and $W_{k'}$. Observe that $W_kW_{k'}^{*}=0$ for $|k-k'|\geq 2$. It suffices to estimate $W_k^{*}W_{k'}.$ Assume $k\geq k'$. By Remark \ref{Rmk polynomial type func}, it is easily verified that $|D(x,y)|^z$ is, uniformly in $y$, a polynomial type function with respect to $x$ on the interval $|x|\leq 2^{(k-N_0)/a_1}$. By Lemma
\ref{lemma almost orthogonal esti}, we obtain
\begin{align*}
\|W_{k}^{*}W_{k'}\|&\leq C\left[|\lambda|2^{k\cdot\sum\limits_{l=1}^n\textbf{N}\begin{bmatrix}
\begin{smallmatrix}
\cdot\\
l
\end{smallmatrix}
\end{bmatrix}}\right]^{-1}\left(2^{k\cdot\sum\limits_{l=1}^n\textbf{N}\begin{bmatrix}
\begin{smallmatrix}
\cdot\\
l
\end{smallmatrix}
\end{bmatrix}}\right)^{1/2}\left(2^{k'\cdot\sum\limits_{l=1}^n\textbf{N}\begin{bmatrix}
\begin{smallmatrix}
\cdot\\
l
\end{smallmatrix}
\end{bmatrix}}\right)^{1/2}
\leq C|\lambda|^{-1}2^{-|k-k'|\sum\limits_{l=1}^n\textbf{N}\begin{bmatrix}
\begin{smallmatrix}
\cdot\\
l
\end{smallmatrix}
\end{bmatrix}}.
\end{align*}
This implies the desired estimate (\ref{Section 3 main L2 decay estimate}).\\

\textbf{Case (ii)}\quad $k\ll a_nj.$\\

On the support of $W_{j,k}$, we have
\begin{align*}
&|D(x,y)|\approx 2^{jm}2^{j\sum\limits_{l=1}^r\textbf{N}\begin{bmatrix}
\begin{smallmatrix}
\cdot\\
l
\end{smallmatrix}
\end{bmatrix}a_l}=2^{jA_r},~~~~~\left|S_{xy}''(x,y)\right|\approx 2^{jm}2^{ks}2^{j\sum\limits_{l=1}^n\textbf{N}\begin{bmatrix}
\begin{smallmatrix}
\cdot\\
l
\end{smallmatrix}
\end{bmatrix}a_l}.
\end{align*}
Assume first $s>0$, $j\geq j'$ and $|k-k'|\leq 1$. The oscillation estimate is
\[
\|W_{j,k}W_{j',k'}^{*}\|\leq C\left[|\lambda|2^{jm}2^{ks}2^{j\sum\limits_{l=1}^n\textbf{N}\begin{bmatrix}
\begin{smallmatrix}
\cdot\\
l
\end{smallmatrix}
\end{bmatrix}a_l}\right]^{-1}2^{jA_r\Re(z)}2^{j'A_r\Re(z)}.
\]
The size estimate is
\[
\|W_{j,k}W_{j',k'}^{*}\|\leq C2^{jA_r\Re(z)}2^{j'A_r\Re(z)}2^{j/2}2^{k/2}2^{j'/2}2^{k'/2}.
\]
Taking a convex combination of the above two estimates yields the following terms and their corresponding exponents:
\begin{align*}
&2^j:\quad -\left(A_r+\sum\limits_{l=r+1}^n\textbf{N}\begin{bmatrix}
\cdot\\
l
\end{bmatrix}a_l\right)\theta+A_r\Re(z)+\frac{1}{2}(1-\theta);\\
&2^{j'}:\quad A_r\Re(z)+\frac{1}{2}(1-\theta);\\
&2^k:\quad -\theta s+B_r\theta=\sum\limits_{l=r+1}^n\textbf{N}\begin{bmatrix}
\cdot\\
l
\end{bmatrix}\theta,
\end{align*}
where we identify $k'$ with $k$. By the assumption $k'\ll a_nj$, we have
\[
2^{k\cdot\sum\limits_{l=r+1}^n\textbf{N}\begin{bmatrix}
\begin{smallmatrix}
\cdot\\
l
\end{smallmatrix}
\end{bmatrix}\theta}\lesssim 2^{j'\cdot\sum\limits_{l=r+1}^n\textbf{N}\begin{bmatrix}
\begin{smallmatrix}
\cdot\\
l
\end{smallmatrix}
\end{bmatrix}a_n\theta}.
\]
Then the resulting exponent of $2^{j'}$ becomes
\[
A_r\Re(z)+\sum\limits_{l=r+1}^n\textbf{N}\begin{bmatrix}
\cdot\\
l
\end{bmatrix}a_n\theta+\frac{1}{2}(1-\theta).
\]
To prove the desired estimate
\[
\|W_{j,k}W_{j',k'}^{*}\|\leq C|\lambda|^{-\frac{1}{1+B_r}}2^{-|j-j'|\delta}
\]
for some $\delta>0$, it suffices to show that the sum of the exponents of $2^j$ and $2^{j'}$ is nonnegative. Indeed, it is easily verified that
\begin{align*}
&-\left(A_r+\sum\limits_{l=r+1}^n\textbf{N}\begin{bmatrix}
\cdot\\
l
\end{bmatrix}a_l\right)\theta+2A_r\Re(z)+\sum\limits_{l=r+1}^n\textbf{N}\begin{bmatrix}
\cdot\\
l
\end{bmatrix}a_n\theta+(1-\theta)\\
&=-A_r\theta-\sum\limits_{l=r+1}^n\textbf{N}\begin{bmatrix}
\cdot\\
l
\end{bmatrix}a_l\theta+\left(A_r-B_r\right)\theta
+\sum\limits_{l=r+1}^n\textbf{N}\begin{bmatrix}
\cdot\\
l
\end{bmatrix}a_n\theta+B_r\theta\\
&=\sum\limits_{l=r+1}^n\textbf{N}\begin{bmatrix}
\cdot\\
l
\end{bmatrix}(a_n-a_l)\theta\geq 0.
\end{align*}
Hence we have established our claim. By the same argument as above, we can prove that $W_{j,k}^{\ast}W_{j',k'}$ satisfies a similar estimate.

Consider the case $s=0$. Now define the operator $W_j$ by $W_j=\sum_kW_{j,k}$
with the summation taken over $k\ll a_nj$ for fixed $j$. Similarly, the oscillation estimate and the size estimate are given by
\[\|W_{j}W_{j'}^{*}\|\leq C\left[|\lambda|2^{jm}2^{j\sum\limits_{l=1}^n\textbf{N}\begin{bmatrix}
\begin{smallmatrix}
\cdot\\
l
\end{smallmatrix}
\end{bmatrix}a_l}\right]^{-1}2^{jA_r\Re(z)}2^{j'A_r\Re(z)}, \quad j\geq j',\]
and
\[
\|W_{j}W_{j'}^{*}\|\leq C2^{jA_r\Re(z)}2^{j'A_r\Re(z)}2^{j/2}2^{a_nj/2}2^{j'/2}2^{a_nj'/2}.
\]
A convex combination gives $\|W_{j}W_{j'}^{*}\|\leq C|\lambda|^{-\theta}2^{-|j-j'|\delta}$ with some $\delta>0$. Indeed, the corresponding terms and their exponents here are given by
\begin{align*}
&2^j:\quad -\left(A_r+\sum\limits_{l=r+1}^n\textbf{N}\begin{bmatrix}
\cdot\\
l
\end{bmatrix}a_l\right)\theta+A_r\Re(z)+\frac{1+a_n}{2}(1-\theta);\\
&2^{j'}:\quad A_r\Re(z)+\frac{1+a_n}{2}(1-\theta).
\end{align*}
The sum of these two exponents equals
\begin{align*}
&-\left(A_r+\sum\limits_{l=r+1}^n\textbf{N}\begin{bmatrix}
\cdot\\
l
\end{bmatrix}a_l\right)\theta+2A_r\Re(z)+(1+a_n)(1-\theta)\\
&=-\left(A_r+\sum\limits_{l=r+1}^n\textbf{N}\begin{bmatrix}
\cdot\\
l
\end{bmatrix}a_l\right)\theta+(A_r-B_r)\theta+(1+a_n)B_r\theta\\
&=\sum\limits_{l=r+1}^n\textbf{N}\begin{bmatrix}
\cdot\\
l
\end{bmatrix}(a_n-a_l)\theta\geq 0.
\end{align*}
Combining all above results, we have obtained the desired estimate in Case (ii).\\

\textbf{Case (iii)}\quad $k= a_{t_1}j+O(1)$, i.e., $|k-a_{t_1}j|\leq N_0$ for some $t_1\in \{1,2,\cdots,n\}.$\\

In what follows, our proof will be divided into two subcases $t_1>r$ and $1\leq t_1\leq r$. \\

\textbf{Subcase (iii--1) $ t_1 > r$.}

On the support of $W_{j,k}$,
$
|D(x,y)|\approx 2^{jm}2^{j\sum\limits_{l=1}^r\textbf{N}\begin{bmatrix}
\begin{smallmatrix}
\cdot\\
l
\end{smallmatrix}
\end{bmatrix}a_l}=2^{jA_r}.
$
By the decomposition procedure described as above, we shall obtain
$
W_{j,k}=\sum W_{j,k,\mu_1,\cdots,\mu_N}^{\sigma_1,\sigma_2,\cdots,\sigma_N},
$
where each operator in this summation is nondegenerate.

Assume $\mu_N$ lies in the range $a_{t_1,t_2,\cdots,t_{N+1}+1}^{\alpha_1,\cdots,\alpha_N}j\ll \mu_N\ll a_{t_1,\cdots,t_{N+1}}^{\alpha_1,\cdots,\alpha_N}j.$
We can see that the size of $S_{xy}''$ is comparable to
\begin{align*}
\left|S_{xy}''(x,y)\right|\approx& 2^{jm}2^{ks}\prod_{a=1}^N\left(\prod_{t< t_a}2^{j\textbf{N}\begin{bmatrix}
\begin{smallmatrix}
\alpha_1,\cdots,\alpha_{a-1},\cdot\\
t_1,\cdots,t_{a-1},t
\end{smallmatrix}
\end{bmatrix}a_{t_1,\cdots,t_{a-1},t}^{\alpha_1,\cdots,\alpha_{a-1}}}
\right)\cdot\left(\prod_{t>t_a}
2^{\mu_{a-1}\textbf{N}\begin{bmatrix}
\begin{smallmatrix}
\alpha_1,\cdots,\alpha_{a-1},\cdot\\
t_1,\cdots,t_{a-1},t
\end{smallmatrix}
\end{bmatrix}}\right)\cdot\\
&\left(\prod_{\alpha< \alpha_a}2^{\mu_{a-1}\textbf{N}\begin{bmatrix}
\begin{smallmatrix}
\alpha_1,\cdots,\alpha_{a-1},\alpha\\
t_1,\cdots,t_{a-1},t_a
\end{smallmatrix}
\end{bmatrix}}\right)\cdot
\left(\prod_{\alpha> \alpha_a}2^{j\textbf{N}\begin{bmatrix}
\begin{smallmatrix}
\alpha_1,\cdots,\alpha_{a-1},\alpha\\
t_1,\cdots,t_{a-1},t_a
\end{smallmatrix}
\end{bmatrix}a_{t_1,\cdots,t_{a-1},t_a}^{\alpha_1,\cdots,\alpha_{a-1}}}\right)\cdot\\
&\left(\prod_{t\leq t_{N+1}}2^{j\textbf{N}\begin{bmatrix}
\begin{smallmatrix}
\alpha_1,\cdots,\alpha_{N},\cdot\\
t_1,\cdots,t_{N},t
\end{smallmatrix}
\end{bmatrix}a_{t_1,\cdots,t_{N},t}^{\alpha_1,\cdots,
\alpha_{N}}}\right)
\cdot
\left(\prod_{t> t_{N+1}}2^{\mu_{N}\textbf{N}
\begin{bmatrix}
\begin{smallmatrix}
\alpha_1,\cdots,\alpha_{N},\cdot\\
t_1,\cdots,t_{N},t
\end{smallmatrix}
\end{bmatrix}}\right).
\end{align*}
By the operator van der Corput Lemma \ref{Lemma operator vers van der Corput}, we get
\[
\left\|W_{j,k,\mu_1,\cdots,\mu_N}^{\sigma_1,\sigma_2,\cdots,\sigma_N}\right\|\leq C\Big(|\lambda|\inf\left|S_{xy}''(x,y)\right|\Big)^{-1/2}\left(2^{jA_r}\right)^{\Re(z)}
\]
with the infimum taken over the support of $W_{j,k,\mu_1,\cdots,\mu_N}^{\sigma_1,\sigma_2,\cdots,\sigma_N}$.
The size estimate gives
\[
\left\|W_{j,k,\mu_1,\cdots,\mu_N}^{\sigma_1,\sigma_2,\cdots,\sigma_N}\right\|\leq
\left(2^{jA_r}\right)^{\Re(z)}2^{\mu_N/2}\left(2^{\mu_N}2^{-j(a_{t_1}-1)}\right)^{1/2}.
\]

Let $\theta=(1+B_r)^{-1}$. If we take a convex combination of the above two estimates, we see that the exponent of $2^{\mu_N}$ is
\begin{eqnarray*}
-\frac12\sum_{t\geq t_{N+1}}\begin{bmatrix}
\alpha_1,\cdots,\alpha_N,\cdot\\
t_1,\cdots,t_N,t
\end{bmatrix}\theta + 1-\theta
\geq
-\frac12 B_r\theta + B_r\theta \geq 0.
\end{eqnarray*}
Since $a_{t_1,t_2,\cdots,t_{N+1}+1}^{\alpha_1,\cdots,\alpha_N}j\ll \mu_N\ll a_{t_1,\cdots,t_{N+1}}^{\alpha_1,\cdots,\alpha_N}j\leq 0$, we can apply Lemma \ref{section 2 control of abs sum} to the summation $W_{j,k}=\sum_{\mu_N}W_{j,k,\mu_1,\cdots,\mu_N}^{\sigma_1,\cdots,
\sigma_N}$. Therefore it suffices to show that each term of the convex combination is bounded by a constant multiple of $|\lambda|^{ -\frac{1}{ 2(1+B_r)} }$. In fact, by the fact $\mu_N\leq \mu_1\wedge\mu_2\wedge\cdots\wedge\mu_{N-1}+O(1)$, the $L^2$ operator norm $\|W_{j,k,\mu_1,\cdots,\mu_N}^{\sigma_1,\cdots,
\sigma_N}\|$ is bounded by a constant multiple of
\begin{eqnarray*}
&&\left[\Big(|\lambda|\inf |S_{xy}''(x,y)|\Big)^{-1/2}
(2^{jA_r})^{\Re(z)}\right]^{\theta}
\left[(2^{jA_r})^{\Re(z)} 2^{\mu_N/2}
\left(2^{\mu_N} 2^{-j(a_{t_1}-1)}  \right)^{1/2}\right]^{1-\theta}\\
&\leq&
C \left( |\lambda| 2^{\mu_N(A_r/a_{t_1\cdots t_{N+1}}^{\alpha_1\cdots\alpha_N}+B_r)}
\right)^{-\theta/2}
\left(2^{\mu_N A_r/a_{t_1\cdots t_{N+1}}^{\alpha_1\cdots\alpha_N}}
\right)^{\Re(z)}
\left(2^{\mu_N/2} 2^{\mu_N A_r/2a_{t_1\cdots t_{N+1}}^{\alpha_1\cdots\alpha_N}}
\right)^{1-\theta}\\
&\leq &
C |\lambda|^{ -\frac{1}{ 2(1+B_r)} }
\end{eqnarray*}
since the exponent of $2^{\mu_N}$ in the last inequality is equal to
\begin{equation*}
-\frac{A_r}{a_{t_1\cdots t_{N+1}}^{\alpha_1\cdots\alpha_N}}
\cdot\frac{\theta}{2}
-\frac{B_r\theta}{2}
+\frac{1}{a_{t_1\cdots t_{N+1}}^{\alpha_1\cdots\alpha_N}}
(A_r-B_r)\cdot \frac{\theta}{2}
+\frac{B_r\theta}{2}
+\frac{B_r\theta}{2a_{t_1\cdots t_{N+1}}^{\alpha_1\cdots\alpha_N}}
=0.
\end{equation*}
The above argument applies without essential change for other ranges of $\mu_N$.\\

\textbf{Subcase (iii--2) $1\leq t_1\leq r$.}

The argument is more delicate in this subcase. In what follows, our treatment is somewhat different depending on whether $a_{t_1}\geq 1$ or not.

Assume now $a_{t_1}\geq 1$. We first insert $\Phi(\sigma_1(y-C_{t_1}^{\alpha_1}x^{a_{t_1}})/2^{\mu_1})$ into $W_{j,k}$ and then obtain smaller pieces $W_{j,k,\mu_1}^{\sigma_1}$. If $\mu_1$ satisfies $\mu_1\geq k-N_0$ for some large positive integer $N_0$, then we choose an arbitrary $\alpha_1'\neq \alpha_1$ and further decompose $W_{j,k,\mu_1}^{\sigma_1}$ into smaller pieces $W_{j,k,\mu_1,\mu_2}^{\sigma_1,\sigma_2}$ by insertion of $\Phi(\sigma_2(y-C_{t_1}^{\alpha_1'}x^{a_{t_1}})/2^{\mu_2})$ into $W_{j,k,\mu_1}^{\sigma_1}$. Without loss of generality, we only consider $\mu_1\ll k$ such that on the support of $W_{j,k,\mu_1}^{\sigma_1}$, $\left|y-C_{t_1}^{\alpha}x^{a_{t_1}}\right|\approx 2^k$ for $\alpha\neq \alpha_1$. Note that $\mu_1\ll k$ also implies $C_{t_1}^{\alpha_1}\in\bR$ if $N_0$ is sufficiently large. Now there are four different ranges for $\mu_1$:
\begin{enumerate}
\item[{\rm (i)}] $a_{t_1,1}^{\alpha_1}j\ll \mu_1\ll a_{t_1}j$;
\item[{\rm (ii)}] $a_{t_1,t+1}^{\alpha_1}j\ll \mu_1\ll a_{t_1,t}^{\alpha_1}j,\quad 1\leq t\leq \textbf{N}\begin{bmatrix}
    \begin{smallmatrix}
\alpha_1\\
t_1
\end{smallmatrix}
\end{bmatrix}-1$;
\item[{\rm (iii)}] $\mu_1\ll a_{t_1,t}^{\alpha_1}j$ with $t=\textbf{N}\begin{bmatrix}
    \begin{smallmatrix}
\alpha_1\\
t_1
\end{smallmatrix}
\end{bmatrix}$;
\item[{\rm (iv)}] $\mu_1=a_{t_1,t}^{\alpha_1}j + O(1)$
for some $1\leq t\leq \textbf{N}\begin{bmatrix}
\begin{smallmatrix}
\alpha_1\\
t_1
\end{smallmatrix}
\end{bmatrix}$.
\end{enumerate}
In the previous three cases, each operator $W_{j,k,\mu_1}^{\sigma_1}$ is nondegenerate and our resolution of singularities is complete. Hence we need only consider the case (iv). More precisely, we assume
$\left|\mu_1-a_{t_1,t_2}^{\alpha_1}j\right|\lesssim 1$. Then we further insert $\Phi(\sigma_2(y-C_{t_1}^{\alpha_1}x^{a_{t_1}}
-C_{t_1,t_2}^{\alpha_1,\alpha_2}x^{a_{t_1,t_2}^{\alpha_1}})/{2^{\mu_2}})$
into $W_{j,k,\mu_1}^{\sigma_1}$ and denote this new operator by $W_{j,k,\mu_1,\mu_2}^{\sigma_1,\sigma_2}$. For the same reason as for $\mu_1$, we may assume $\mu_2\ll \mu_1$
such that for $\alpha\neq \alpha_2$,
$
|y-C_{t_1}^{\alpha_1}x^{a_{t_1}}
-C_{t_1,t_2}^{\alpha_1,\alpha}x^{a_{t_1,t_2}^{\alpha_1}}|\approx 2^{\mu_1}
$
in the support of $W_{j,k,\mu_1,\mu_2}^{\sigma_1,\sigma_2}$.\\

Generally, if $W_{j,k,\mu_1,\cdots,\mu_p}^{\sigma_1,\sigma_2,\cdots,\sigma_p}$ has been defined and $\mu_p$ lies in the range
$\mu_p=a_{t_1,\cdots,t_{p+1}}^{\alpha_1,\cdots,\alpha_p}j+O(1),$
we choose an index $\alpha_{p+1}$ and insert
$
\Phi\Big(\sigma_{p+1}(y
-\sum_{i=1}^{p+1}C_{t_1,\cdots,t_{i}}^{\alpha_1,\cdots,\alpha_i}
x^{a_{t_1,\cdots,t_i}^{\alpha_1,\cdots,\alpha_{i-1}}})/{2^{\mu_{p+1}}}\Big)
$
into $W_{j,k,\mu_1,\cdots,\mu_p}^{\sigma_1,\cdots,\sigma_p}$. In this way, we obtain a new operator $W_{j,k,\mu_1,\cdots,\mu_{p+1}}^{\sigma_1,\cdots,\sigma_{p+1}}$. Similar to our discussion for $\mu_1$, there are also four possible ranges for $\mu_{p+1}$.

Now we assume the decomposition process stops in the $N$-th step. Then $\mu_N$ must satisfy one of the following restrictions:
\begin{enumerate}
\item[{\rm (i)}] $a_{t_1,\cdots,t_N,1}^{\alpha_1,\cdots,\alpha_N}j\ll \mu_N\ll \mu_{N-1}j$;
\item[{\rm (ii)}] $a_{t_1,\cdots,t_N,t_{N+1}+1}^{\alpha_1,\cdots,\alpha_N}j\ll \mu_N\ll a_{t_1,\cdots,t_{N+1}}^{\alpha_1,\cdots,\alpha_N}j$ for some
    $1\leq t_{N+1}\leq \textbf{N}\begin{bmatrix}
\begin{smallmatrix}
\alpha_1,\cdots,\alpha_N\\
t_1,\cdots,t_N
\end{smallmatrix}
\end{bmatrix}-1$ and\\
$1\leq \alpha_N\leq \textbf{N}\begin{bmatrix}
\begin{smallmatrix}
\alpha_1,\cdots,\alpha_{N-1}, \cdot\\
t_1,\cdots,t_{N-1}, t_N
\end{smallmatrix}
\end{bmatrix}$;
\item[{\rm (iii)}] $\mu_N\ll a_{t_1,\cdots,t_{N+1}}^{\alpha_1,\cdots,\alpha_N}j$ with $t_{N+1}=\textbf{N}\begin{bmatrix}
\begin{smallmatrix}
\alpha_1,\cdots,\alpha_N\\
t_1,\cdots,t_N
\end{smallmatrix}
\end{bmatrix}$;
\item[{\rm (iv)}] For some $t$ we have $\mu_N= a_{t_1,\cdots,t_N,t}^{\alpha_1,\cdots,\alpha_N}j + O(1)$, and for all $\alpha$
    \begin{equation*}
    \Big|y-\sum_{i=1}^N C_{t_1,\cdots,t_i}^{\alpha_1,\cdots,\alpha_{i}} x^{a_{t_1\cdots t_i}^{\alpha_1\cdots\alpha_{i-1}}}
    -
     C_{t_1,\cdots,t_N,t}^{\alpha_1,\cdots,\alpha_N,\alpha} x^{a_{t_1\cdots t_Nt}^{\alpha_1\cdots\alpha_N}}\Big|
     \approx
     2^{\mu_N},~~~(x,y)\in \supp(W_{j,k,\mu_1,\cdots,\mu_N}^{\sigma_1,\cdots,\sigma_N}).
    \end{equation*}
\end{enumerate}
By our assumption, $W_{j,k,\mu_1,\cdots,\mu_N}^{\sigma_1,\cdots,\sigma_N}$ is nondegenerate such that both the Hessian $S_{xy}''$ and the damping factor $D$ are bounded from below and above by positive constants.

Now we are going to prove the desired $L^2$ decay estimate for $\sum W_{j,k,\mu_1,\cdots,\mu_N}^{\sigma_1,\cdots,\sigma_N}$ with the summation taken over the above range (ii). Other cases can be treated in the same way.
Before our application of the almost orthogonality estimate in Lemma \ref{lemma almost orthogonal esti}, we shall verify all assumptions there for $ W_{j,k,\mu_1,\cdots,\mu_N}^{\sigma_1,\cdots,\sigma_N}$ and $ W_{j,k,\mu_1,\cdots,\mu_N'}^{\sigma_1,\cdots,\sigma_N}$. Consider the domain
$\Omega_{j,k,\mu_1,\cdots,\mu_N}^{\sigma_1,\cdots,\sigma_N}$ which consists of all points $(x,y)$ satisfying
\begin{eqnarray*}
2^{j-1}\leq x \leq 2^{j+1}, \quad 2^{k-1}\leq y \leq 2^{k+1},
\quad 2^{\mu_t-1}\leq \sigma_t \left(y-\sum_{i=1}^{t} C_{t_1,\cdots,t_i}^{\alpha_1,\cdots,\alpha_{i}}x^{a_{t_1\cdots t_i}^{\alpha_1\cdots\alpha_{i-1}}}\right)
\leq 2^{\mu_t+1}
\end{eqnarray*}
for $1\leq t \leq N$. It is easy to see that this domain is horizontally and vertically convex. We also define its expanded domain by
\begin{eqnarray*}
& &\Omega_{j,k,\mu_1,\cdots,\mu_N}^{\sigma_1,\cdots,\sigma_N~\ast}:
~2^{j-1}-\epsilon 2^{j} \leq x \leq 2^{j+1}+\epsilon 2^{j}, \quad 2^{k-1}-\epsilon 2^{k}\leq y \leq 2^{k+1}+\epsilon 2^{k}, \\
& &2^{\mu_t-1}-\epsilon 2^{\mu_t}\leq \sigma_t \left(y-\sum_{i=1}^{t} C_{t_1,\cdots,t_i}^{\alpha_1,\cdots,\alpha_{i}}x^{a_{t_1\cdots t_i}^{\alpha_1\cdots\alpha_{i-1}}}\right)
\leq 2^{\mu_t+1}+\epsilon 2^{\mu_t}~~{\rm for}~~1\leq t \leq N
\end{eqnarray*}
for sufficiently small $\epsilon>0$. By direct verification as in (\ref{section 2 nonempty intersection property I}) and (\ref{section 2 nonempty intersection property II}), one can see that $\Omega_{j,k,\mu_1,\cdots,\mu_N}^{\sigma_1,\cdots,\sigma_N~\ast}$
satisfies all assumptions in Definitions \ref{Def expanded domains I} and \ref{Def expanded domains II} and hence it becomes an expanded domain for $\Omega_{j,k,\mu_1,\cdots,\mu_N}^{\sigma_1,\cdots,\sigma_N}$. For sufficiently small $\epsilon$, all assumptions in Lemma \ref{lemma almost orthogonal esti} are also true for the above expanded domain.

By Lemma \ref{lemma almost orthogonal esti}, the $L^2$ operator norm of $W_{j,k,\mu_1,\cdots,\mu_N}^{\sigma_1,\cdots,\sigma_N}
W_{j,k,\mu_1,\cdots,\mu_N^{'}}^{\sigma_1,\cdots,\sigma_N~*}$ is bounded by a constant multiple of
\begin{align}\label{Section 3 general almost ortho osci esti}
&\Bigg[|\lambda|2^{jA_{t_1-1}}2^{kB_{t_1}}
\prod_{\alpha\neq \alpha_1}2^{k\textbf{N}\begin{bmatrix}
\begin{smallmatrix}
\alpha\\
t_1
\end{smallmatrix}
\end{bmatrix}}
\prod_{l<t_2}2^{\textbf{N}\begin{bmatrix}
\begin{smallmatrix}
\alpha_1,\cdot\\
t_1,l
\end{smallmatrix}
\end{bmatrix}a_{t_1,l}^{\alpha_1}j}
\prod_{l>t_2}2^{\textbf{N}\begin{bmatrix}
\begin{smallmatrix}
\alpha_1,\cdot\\
t_1,l
\end{smallmatrix}
\end{bmatrix}\mu_1}
\prod_{\alpha\neq\alpha_2}2^{\mu_1\textbf{N}\begin{bmatrix}
\begin{smallmatrix}
\alpha_1,\alpha\\
t_1,t_2
\end{smallmatrix}
\end{bmatrix}}\cdots  \nonumber  \\
&\prod_{l\leq t_{N+1}}2^{j\textbf{N}\begin{bmatrix}
\begin{smallmatrix}
\alpha_1,\cdots,\alpha_N,\cdot\\
t_1,\cdots,t_N,l
\end{smallmatrix}
\end{bmatrix}a_{t_1,\cdots,t_N,l}^{\alpha_1,\cdots,\alpha_N}}\prod_{l> t_{N+1}}2^{\textbf{N}\begin{bmatrix}
\begin{smallmatrix}
\alpha_1,\cdots,\alpha_N,\cdot\\
t_1,\cdots,t_N,l
\end{smallmatrix}
\end{bmatrix}\mu_N}\Bigg]^{-1}
\Bigg[2^{jm}2^{\sum\limits_{l<t_1}\textbf{N}\begin{bmatrix}
\begin{smallmatrix}
\cdot\\
l
\end{smallmatrix}
\end{bmatrix}a_lj}\cdot   \nonumber  \\
&2^{k\sum\limits_{r\geq l>t_1}\textbf{N}\begin{bmatrix}
\begin{smallmatrix}
\cdot\\
l
\end{smallmatrix}
\end{bmatrix}}\cdot\prod_{\alpha\neq \alpha_1}2^{k\textbf{N}\begin{bmatrix}
\begin{smallmatrix}
\alpha\\
t_1
\end{smallmatrix}
\end{bmatrix}}
\prod_{l<t_2}2^{\textbf{N}\begin{bmatrix}
\begin{smallmatrix}
\alpha_1,\cdot\\
t_1,l
\end{smallmatrix}
\end{bmatrix}a_{t_1,l}^{\alpha_1}j}2^{\sum\limits_{l>t_2}
\textbf{N}\begin{bmatrix}
\begin{smallmatrix}
\alpha_1,\cdot\\
t_1,l
\end{smallmatrix}
\end{bmatrix}\mu_1}\cdot\prod_{\alpha\neq\alpha_2}2^{\mu_1\textbf{N}
\begin{bmatrix}
\begin{smallmatrix}
\alpha_1,\alpha\\
t_1,t_2
\end{smallmatrix}
\end{bmatrix}}\cdot\cdots\cdot    \nonumber   \\
&\prod_{l\leq t_{N+1}}2^{j\textbf{N}\begin{bmatrix}
\begin{smallmatrix}
\alpha_1,\cdots,\alpha_N,\cdot\\
t_1,\cdots,t_N,l
\end{smallmatrix}
\end{bmatrix}a_{t_1,\cdots,t_N,l}^{\alpha_1,\cdots,\alpha_N}}
\prod_{l> t_{N+1}}2^{\textbf{N}\begin{bmatrix}
\begin{smallmatrix}
\alpha_1,\cdots,\alpha_N,\cdot\\
t_1,\cdots,t_N,l
\end{smallmatrix}
\end{bmatrix}\frac{\mu_N+\mu_{N}'}{2}}\Bigg]^{2\Re(z)}.
\end{align}
The size estimate gives us the following bound of the operator norm of $W_{j,k,\mu_1,\cdots,\mu_N}^{\sigma_1,\cdots,\sigma_N}
W_{j,k,\mu_1,\cdots,\mu_N^{'} }^{\sigma_1,\cdots,\sigma_N~\ast}$.
\begin{align}\label{Section 3 general size estimate}
&C\Bigg[2^{mj}2^{\sum\limits_{l<t_1}\textbf{N}\begin{bmatrix}
\begin{smallmatrix}
\cdot\\
l
\end{smallmatrix}
\end{bmatrix}a_lj}2^{k\sum\limits_{r\geq l>t_1} \textbf{N} \begin{bmatrix}
\begin{smallmatrix}
\cdot\\
l
\end{smallmatrix}
\end{bmatrix}}\cdot\prod_{\alpha\neq \alpha_1}2^{k\textbf{N}\begin{bmatrix}
\begin{smallmatrix}
\alpha\\
t_1
\end{smallmatrix}
\end{bmatrix}}2^{\sum\limits_{l<t_2}\textbf{N}\begin{bmatrix}
\begin{smallmatrix}
\alpha_1,\cdot\\
t_1,l
\end{smallmatrix}
\end{bmatrix}a_{t_1,l}^{\alpha_1}j}2^{\sum\limits_{l>t_2}\textbf{N}\begin{bmatrix}
\begin{smallmatrix}
\alpha_1,\cdot\\
t_1,l
\end{smallmatrix}
\end{bmatrix}\mu_1}\cdot  \nonumber \\
&\prod_{\alpha\neq\alpha_2}2^{\mu_1\textbf{N}\begin{bmatrix}
\begin{smallmatrix}
\alpha_1,\alpha\\
t_1,t_2
\end{smallmatrix}
\end{bmatrix}}\cdot\cdots\cdot \prod_{l\leq t_N}2^{j\textbf{N}\begin{bmatrix}
\begin{smallmatrix}
\alpha_1,\cdots,\alpha_{N-1},\cdot\\
t_1,\cdots,t_{N-1},l
\end{smallmatrix}
\end{bmatrix}a_{t_1,\cdots,t_{N-1},l}^{\alpha_1,\cdots,\alpha_{N-1}}}
2^{\mu_{N-1}\sum\limits_{l>t_N}\textbf{N}\begin{bmatrix}
\begin{smallmatrix}
\alpha_1,\cdots,\alpha_{N-1},\cdot\\
t_1,\cdots,t_{N-1},l
\end{smallmatrix}
\end{bmatrix}}\nonumber \\
&\prod_{\alpha\neq \alpha_N}2^{\mu_{N-1}\textbf{N}\begin{bmatrix}
\begin{smallmatrix}
\alpha_1,\cdots,\alpha_{N-1},\cdot\\
t_1,\cdots,t_{N-1},t_N
\end{smallmatrix}
\end{bmatrix}}
2^{j\sum\limits_{l\leq t_{N+1}}\textbf{N}\begin{bmatrix}
\begin{smallmatrix}
\alpha_1,\cdots,\alpha_{N},\cdot\\
t_1,\cdots,t_{N},l
\end{smallmatrix}
\end{bmatrix} a_{t_1,\cdots, t_N,l}^{\alpha_1,\cdots,\alpha_N} }2^{\sum\limits_{l\geq t_{N+1}+1}\textbf{N}\begin{bmatrix}
\begin{smallmatrix}
\alpha_1,\cdots,\alpha_{N},\cdot\\
t_1,\cdots,t_{N},l
\end{smallmatrix}
\end{bmatrix}\frac{\mu_N+\mu_N'}{2}}
\Bigg]^{2\Re(z)} \nonumber \\
&2^{\mu_N/2}\left(2^{\mu_N}2^{-j(a_{t_1}-1)}\right)^{1/2}\cdot 2^{\mu_N'/2}\left(2^{\mu_N'}2^{-j(a_{t_1}-1)}\right)^{1/2}.
\end{align}
By a convex combination, we see that the $L^2$ operator norm of
$W_{j,k,\mu_1,\cdots,\mu_N}^{\sigma_1,\cdots,\sigma_N}
W_{j,k,\mu_1,\cdots,\mu_N'}^{\sigma_1,\cdots,\sigma_N~\ast}$ is bounded by a constant multiple of the product of terms $2^{-j\delta_j}$, $2^{-k\delta_k}$, $2^{-\mu_1\delta_{\mu_1}}$, $\cdots$, $2^{-\mu_N\delta_{\mu_N}}$ and $2^{-\mu_{N'}\delta_{\mu_{N'}}}$. Since $a_{t_1}\geq 1$, it is clear that $\delta_j\geq 0$. Then it follows that $2^{-j\delta_j}$ is bounded by a constant multiple of
\begin{align}\label{section 3 estimate of 2^j}
& \Bigg(2^{-k\frac{A_{t_1-1}}{a_{t_1}} }2^{-\mu_1\sum\limits_{l<t_2}
\textbf{N}\begin{bmatrix}
\begin{smallmatrix}
\alpha_1,\cdot\\
t_1,l
\end{smallmatrix}
\end{bmatrix}\frac{a_{t_1,l}^{\alpha_1}}{a_{t_1,t_2}^{\alpha_1}}}
\cdots
2^{-\mu_{N}\sum\limits_{t\leq t_{N+1}}\textbf{N}\begin{bmatrix}
\begin{smallmatrix}
\alpha_1,\cdots,\alpha_{N},\cdot\\
t_1,\cdots,t_{N},t
\end{smallmatrix}
\end{bmatrix}\frac{a_{t_1,\cdots,t_{N},t}^{\alpha_1,\cdots,\alpha_{N}}}
{
a_{t_1,\cdots,t_{N},t_{N+1}}^{\alpha_1,\cdots,\alpha_N}}}\Bigg)^{\frac{B_r}{A_r}\theta}2^{-kB_r\theta}
2^{\frac{k}{a_{t_1}}B_r\theta}\nonumber \\
\leq&
C\Bigg(2^{-\mu_N\frac{A_{t_1-1}}{a_{t_1}}-\mu_N\sum\limits_{l<t_2}
\textbf{N}\begin{bmatrix}
\begin{smallmatrix}
\alpha_1,\cdot\\
t_1,l
\end{smallmatrix}
\end{bmatrix}}\cdots 2^{-\mu_{N}\sum\limits_{t\leq t_{N+1}}\textbf{N}\begin{bmatrix}
\begin{smallmatrix}
\alpha_1,\cdots,\alpha_{N-1},\cdot\\
t_1,\cdots,t_{N-1},l
\end{smallmatrix}
\end{bmatrix}}\Bigg)^{\frac{B_r}{A_r}\theta}
2^{-\mu_NB_r\theta}2^{\frac{\mu_N}{a_{t_1}}B_r\theta}.
\end{align}
Similarly, we can give an upper bound for $2^{-k\delta_k}$ as follows.
\begin{align}\label{section 3 estimate of 2^k}
2^{-k\delta_k}&\leq C2^{-\mu_NB_r\theta}
2^{-\mu_N\left(\sum\limits_{t=t_1+1}^r\textbf{N}\begin{bmatrix}
\begin{smallmatrix}
\cdot    \\
l
\end{smallmatrix}
\end{bmatrix}
+\sum\limits_{\alpha\neq\alpha_1}\textbf{N}\begin{bmatrix}
\begin{smallmatrix}
\alpha    \\
t_1
\end{smallmatrix}
\end{bmatrix}\right)(\theta-2\Re(z))}  \nonumber  \\
&\leq C2^{-\mu_NB_r\theta}
2^{-\mu_N\left(\sum\limits_{t=t_1+1}^r\textbf{N}\begin{bmatrix}
\begin{smallmatrix}
\cdot    \\
l
\end{smallmatrix}
\end{bmatrix}+\sum\limits_{\alpha\neq\alpha_1}\textbf{N}\begin{bmatrix}
\begin{smallmatrix}
\alpha      \\
t_1
\end{smallmatrix}
\end{bmatrix}\right)\frac{B_r\theta}{A_r}}.
\end{align}
It is easy to see that $\delta_{\mu_t}\geq 0$ for $1\leq t \leq N-1$. The above decomposition procedure implies
$\mu_p\leq \mu_1\wedge\mu_2\wedge\cdots\wedge\mu_{p-1}+O(1).$ Hence
$2^{-\mu_t\delta_{\mu_t}} \lesssim 2^{-\mu_N\delta_{\mu_t}}.$
Combining this estimate with (\ref{section 3 estimate of 2^j}) and (\ref{section 3 estimate of 2^k}), we see that the sum of the exponents of $2^{\mu_N}$ is equal to
\begin{equation*}
-\Big[\Big( m+\sum\limits_{t<t_1}\textbf{N}\begin{bmatrix}
\begin{smallmatrix}
\cdot    \\
l
\end{smallmatrix}
\end{bmatrix}a_l\Big)\big{/}a_{t_1}
+\sum\limits_{t<t_2}\textbf{N}\begin{bmatrix}
\begin{smallmatrix}
\alpha_1,\cdot      \\
t_1,t
\end{smallmatrix}
\end{bmatrix}
 + \cdots +
\sum\limits_{t\leq t_{N+1}}\textbf{N}\begin{bmatrix}
\begin{smallmatrix}
\alpha_1\cdots\alpha_N\cdot      \\
t_1\cdots t_N t
\end{smallmatrix}
\end{bmatrix}
 \Big]\frac{B_r}{A_r}\theta-
B_r\theta+\frac{1}{a_{t_1}}B_r\theta-
\end{equation*}
\begin{equation*}
B_r\theta-
\Big(
\sum\limits_{t=t_1+1}^r\textbf{N}\begin{bmatrix}
\begin{smallmatrix}
\cdot    \\
l
\end{smallmatrix}
\end{bmatrix}
+
\sum\limits_{\alpha\neq \alpha_1}^r\textbf{N}\begin{bmatrix}
\begin{smallmatrix}
\alpha   \\
t_1
\end{smallmatrix}
\end{bmatrix}
\Big)\frac{B_r}{A_r}\theta-\cdots-
\sum\limits_{t>t_{N+1}}\textbf{N}\begin{bmatrix}
\begin{smallmatrix}
\alpha_1\cdots\alpha_N\cdot   \\
t_1\cdots t_N t
\end{smallmatrix}
\end{bmatrix}\Re(z)
+1-\theta
\end{equation*}
\begin{equation*}
=-\Big[
\Big( m+\sum\limits_{t<t_1}\textbf{N}\begin{bmatrix}
\begin{smallmatrix}
\cdot    \\
l
\end{smallmatrix}
\end{bmatrix}a_l\Big)\big{/}a_{t_1}
+\sum\limits_{t=t_1}^r\textbf{N}\begin{bmatrix}
\begin{smallmatrix}
\cdot      \\
t
\end{smallmatrix}
\end{bmatrix}
\Big]
 \frac{B_r}{A_r}\theta-B_r\theta+\frac{1}{a_{t_1}}B_r\theta
 -\sum\limits_{t>t_{N+1}}\textbf{N}\begin{bmatrix}
\begin{smallmatrix}
\alpha_1\cdots\alpha_N\cdot   \\
t_1\cdots t_N t
\end{smallmatrix}
\end{bmatrix}\Re(z)
\end{equation*}
\begin{equation*}
\geq
-\frac{A_r}{a_{t_1}}\cdot \frac{B_r}{A_r}\theta
-B_r\theta+\frac{1}{a_{t_1}}B_r\theta
 -\sum\limits_{t>t_{N+1}}\textbf{N}\begin{bmatrix}
\begin{smallmatrix}
\alpha_1\cdots\alpha_N\cdot   \\
t_1\cdots t_N t
\end{smallmatrix}
\end{bmatrix}\Re(z)
~~~~~~~~~~~~~~~~~~~~~~~~~~~~~~~~~~~~~~~~~~~~~~
\end{equation*}
\begin{equation*}
=-(1-\theta)-\sum\limits_{t>t_{N+1}}\textbf{N}\begin{bmatrix}
\begin{smallmatrix}
\alpha_1\cdots\alpha_N\cdot   \\
t_1\cdots t_N t
\end{smallmatrix}
\end{bmatrix}\Re(z).
~~~~~~~~~~~~~~~~~~~~~~~~~~~~~~~~~~~~~~~~~~~~~~~~~~~~~~~~
~~~~~~~~~~~~~~
\end{equation*}
As a result, we have
\begin{equation}\label{sec3 general alm ortho}
\left\|W_{j,k,\mu_1,\cdots,\mu_N}^{\sigma_1,\cdots,\sigma_N}
W_{j,k,\mu_1,\cdots,\mu_N'}^{\sigma_1,\cdots,\sigma_N~\ast}\right\|
\leq C|\lambda|^{-\theta}2^{-|\mu_N-\mu_N'|\delta}
\end{equation}
with $\delta$ given by
\begin{equation*}
\delta=
\sum\limits_{t>t_{N+1}}\textbf{N}\begin{bmatrix}
\begin{smallmatrix}
\alpha_1,\cdots,\alpha_N,\cdot   \\
t_1,\cdots, t_N, t
\end{smallmatrix}
\end{bmatrix}\Re(z)
+
(1-\theta)>0.
\end{equation*}
A similar estimate is also true for $W_{j,k,\mu_1,\cdots,\mu_N}^{\sigma_1,\cdots,\sigma_N~\ast}
W_{j,k,\mu_1,\cdots,\mu_N'}^{\sigma_1,\cdots,\sigma_N}.$

It remains to treat the case $0<a_{t_1}<1$. The above argument, with $a_{t_1}\geq1$, does not apply here and some modifications are needed. In fact, the above size estimate produces a term $2^{-j(a_{t_1}-1)}$ and the exponent of $2^j$, appearing in the convex combination of the oscillation estimate (\ref{Section 3 general almost ortho osci esti}) and the size estimate (\ref{Section 3 general size estimate}), may not be negative for $0<a_{t_1}<1$. Hence the estimate (\ref{section 3 estimate of 2^j}) is not generally true.  Since $S$ is a real-analytic function near the origin, there exists a non-vanishing real-analytic function $V$ such that
\begin{equation}\label{sec3 conjuate invariant product}
\prod_{r_{\nu}\in \Psi\begin{bmatrix}
\begin{smallmatrix}
\cdot\\
t_1
\end{smallmatrix}
\end{bmatrix}}(y-r_\nu(x))=V(x,y)\prod(x-h(y)),
\end{equation}
where $h(y)$ is a Puiseux series in $y$ of form
\[
h(y)=d_{l_1}^{\alpha_1}y^{b_{l_1}}
+d_{l_1,l_2}^{\alpha_1,\alpha_2}y^{b_{l_1,l_2}^{\alpha_1}}
+\cdots+d_{l_1,\cdots,l_N}^{\alpha_1,\cdots,\alpha_N}
y^{b_{l_1,\cdots,l_N}^{\alpha_1,\cdots,\alpha_{N-1}}}
+o(y^{b_{l_1,\cdots,l_N}^{\alpha_1,\cdots,\alpha_{N-1}}}),~~~
b_{l_1}=\frac{1}{a_{t_{1}}}.\]
Now we shall prove this equality. By Theorem \ref{Puiseux theorem}, the left side of the equality is conjugation invariant in $x$. Hence the left side is a real-analytic function near the origin; see Lemma 1.2.1 (p18) in \cite{casas-alvero} and Remark \ref{sec3 Remark}. Now it is clear that the equality is an immediate consequence of Theorem \ref{Puiseux theorem} with the roles of $x$ and $y$ changed. Moreover, by the Newton-Puiseux algorithm in the proof of the Puiseux theorem (see \cite{casas-alvero}), we also have $b_{l_1}=\frac{1}{a_{t_{1}}}$.

With the equality (\ref{sec3 conjuate invariant product}), in the case $|k-a_{t_1}j|\lesssim 1$ considered above, we have $|j-k/a_{t_1}|\lesssim 1$. The argument in this case is the same as that of the case $a_{t_1}\geq 1$ with the roles of $x$ and $y$ interchanged. Thus we can also obtain the desired $L^2$ decay estimate. Therefore the proof of Theorem \ref{Thm L2 damping decay estiamte} is complete.
\end{proof}

\section{Damped Oscillatory Integral Operators with Critical Negative Exponents}
Consider the following damped oscillatory integral operator
\begin{equation}\label{damp-opr-1}
W_zf(x)=\int_{-\infty}^{+\infty}e^{i\lambda S(x,y)}|\widetilde{D}(x,y)|^z\varphi(x,y)f(y)dy,
\end{equation}
where $\widetilde{D}$ is slightly different from $D$ in (\ref{section 3 Def damping D}), depending on whether $D$ is of the form
\begin{equation}\label{section 4 translation invar damp fac}
D(x,y)=(y-r_\nu(x))^d
\end{equation}
where $r_\nu(x)=C_1x+o(x)$ for some $C_1\neq 0$.

As in Section 3, we shall decompose $W_z$ into smaller pieces $W_{j,k}$ by insertion of $\Phi(\frac{x}{2^j})\Phi(\frac{y}{2^k})$ into $W_z$. Here, we need only consider $W_z$ in the first quadrant since the same argument applies without essential change in other quadrants.

For clarity, we assume $\varphi$ is supported in the unit square $[-1/2,1/2]^2$. Hence $W_{j,k}=0$ for $j\geq 0$ or $k\geq 0$. Recall that $W_{j,k}$ is defined by
\begin{equation}\label{section 4 Def of Wjk}
W_{j,k}f(x)=\int_{-\infty}^{+\infty}e^{i\lambda S(x,y)}|\widetilde{D}(x,y)|^z
\Phi\l( \frac{x}{2^j} \r)  \Phi\l( \frac{y}{2^k} \r)\varphi(x,y)f(y)dy,~~\Re(z)=-\frac{1}{A_r}.
\end{equation}

Consider first the simplest case $r=0$. Then $D(x,y)=x^m$ and $A_r=m$. In this case, let $\widetilde{D}(x,y)=D(x,y)$. It is clear that $|W_zf(x)|\leq C|x|^{-1}\|f\|_{L^1}.$ Hence $W_z$ is bounded from $L^1$ into $L^{1,\infty}$. As we shall see in Section 5, we can deduce the desired $L^p$ estimate by Lemma
\ref{section 2 stein-weiss interpolation}.

Now assume $r\geq 1$ throughout the rest of this section. In what follows, we will decompose $W_z$ as three parts $W_z^{(1)}$, $W_z^{(2)}$ and $W_z^{(3)}$ as follows.
\begin{equation*}\
W_z^{(1)}=\sum_{k\leq a_rj-N_0} W_{j,k},~~~~~
W_z^{(2)}=\sum_{a_rj-N_0< k < a_1j+N_0} W_{j,k},~~~~~
W_z^{(3)}=\sum_{k\geq a_1j+N_0} W_{j,k},
\end{equation*}
where $N_0$ is a large positive integer. Our choice of damping factors in $W_z^{(i)}$ may be different, depending on whether $D$ is of the form (\ref{section 4 translation invar damp fac}). More precisely, for $W_z^{(1)}$ and $W_z^{(3)}$, we always define
\begin{equation*}
\widetilde{D}(x,y)
=D(x,y)
=x^m \prod_{l=1}^r \prod_{r_{\nu}\in \Psi[\begin{smallmatrix}
\cdot\\
l
\end{smallmatrix}]}(y-r_\nu(x)).
\end{equation*}
For $W_z^{(2)}$, if $D$ is not of the form
(\ref{section 4 translation invar damp fac}) then set $\widetilde{D}(x,y)=D(x,y)$ as above; otherwise, we shall define
\begin{equation}\label{section 4 def damping factor}
\widetilde{D}(x,y)
=\left(|\lambda|2^{kB_1}\right)^{-\frac{\textbf{N}\begin{bmatrix}
\begin{smallmatrix}
\cdot\\
1
\end{smallmatrix}
\end{bmatrix}}{\textbf{N}\begin{bmatrix}
\begin{smallmatrix}
\cdot\\
1
\end{smallmatrix}
\end{bmatrix}+2}}+\prod_{r_{\nu}\in \Psi\begin{bmatrix}
\begin{smallmatrix}
\cdot\\
1
\end{smallmatrix}
\end{bmatrix}}|y-r_\nu(x)|
\end{equation}
and introduce a variant $H^1_E$ of the Hardy space $H^1$ to establish the $H^1_E\rightarrow L^1$ estimate. In the definition of $H^1_E$, we shall need a relevant fact. Recall that $D$ is a real-analytic function, as pointed out in Remark \ref{sec3 Remark}. If $D$ is of the form (\ref{section 4 translation invar damp fac}), then $r_{\nu}$ is a power series with real coefficients since $D$ is real-valued.
\begin{defn}
Let $I_k:=[2^{k-1},2^{k+1}]$ for $k\in\mathbb{Z}$. We say that a measurable function $a$ is an $H_E^1(I_k)$ atom associated with the phase $\lambda S$ if there exists an interval $I\subseteq I_k$ such that
(i) $\supp(a)\subseteq I$; (ii) $\|a\|_{L^\infty}\leq |I|^{-1}$; (iii) $\int_{I}e^{i\lambda S(r_{\nu}^{-1}(c_I),y)}a(y)dy=0$~for some root $r_{\nu}$ in $\Psi\begin{bmatrix}
\begin{smallmatrix}
\cdot\\
1
\end{smallmatrix}
\end{bmatrix}$, where $r_{\nu}^{-1}$ is the inverse of the root $r_{\nu}$ on a small neighborhood of the origin.

The space $H_E^1(I_k)$ is the set of all functions $f\in L^1$, supported in $I_k$, which can be written as
\begin{equation}\label{atom-decom}
f=\sum_{j\in\mathbb{Z}}\lambda_ja_j,
\end{equation}
where each $a_j$ is an $H_E^1(I_k)$ atom. The norm in $H_E^1$ is defined to be
\[\|f\|_{H_E^1}=\inf\Big\{\sum_j|\lambda_j|:f=\sum_j\lambda_ja_j\Big\}\]
with the infimum taken over all sequences $\{\lambda_j\}$ such that \eqref{atom-decom} holds.
\end{defn}

With the above preliminaries, we have the following theorem.
\begin{theorem}\label{Thm Main theorem Sec 4}
Assume $W_z^{(i)}$ are defined as above with $z\in \mathbb{C}$ having real part $\Re(z)=-\frac1{A_r}$. Then there exists a constant $C$ such that
\begin{equation*}
\|W_z^{(1)}f\|_{L^{1,\infty}}\leq C\|f\|_{L^1},~~~~~
\|W_z^{(3)}f\|_{L^1}\leq C\|f\|_{L^1}.
\end{equation*}
If $D$ is not of the form
{\rm (\ref{section 4 translation invar damp fac})}, then $W_z^{(2)}$ is bounded from $L^1$ into itself with the operator norm independent of $\lambda$; otherwise, it is true that
\begin{equation*}
\|W_z^{(2)}f\|_{L^1}
\leq
C\|f\|_{H_E^1(I_k)}.
\end{equation*}
\end{theorem}

\begin{proof}
It is more convenient to divide the argument into several cases.\\

\textrm{$\mathbf{Case~(i)~~k\ll a_rj.}$}

It is clear that $|\widetilde{D}(x,y)|\approx |x|^{A_r}$ for $(x,y)\in \supp(W_{z}^{(1)}).$ Hence $W_z^{(1)}$ is bounded from $L^1$ into $L^{1,\infty}$.\\

\textrm{$\mathbf{Case~(ii)~~k\gg a_1j}$}.

On the support of the operator $W_z^{(3)}$, one has
$|\widetilde{D}(x,y)|\approx x^my^{\sum_{l=1}^r\textbf{N}\begin{bmatrix}
\begin{smallmatrix}
\cdot\\
l
\end{smallmatrix}
\end{bmatrix}}.$
For $|y|\leq 1$,
\begin{align*}
\int_{|y|\gg|x|^{a_1}}|\widetilde{D}(x,y)|^{-\frac{1}{A_r}}dx&\lesssim y^{-\sum\limits_{l=1}^r\textbf{N}\begin{bmatrix}
\begin{smallmatrix}
\cdot\\
l
\end{smallmatrix}
\end{bmatrix}/A_r}y^{\frac{1}{a_1}(1-\frac{m}{A_r})}
\lesssim 1,
\end{align*}
where we have used the fact
$-\sum_{l=1}^r\textbf{N}\begin{bmatrix}
\begin{smallmatrix}
\cdot\\
l
\end{smallmatrix}
\end{bmatrix}/{A_r}+1/{a_1}(1-{m}/{A_r})\geq 0.$
By Fubini's theorem, we see that $W_z^{(3)}$ is bounded on $L^1$.\\

\textrm{$\mathbf{Case~(iii)~~a_{t_0+1}j\ll k\ll a_{t_0}j~~for~some~~1\leq t_0\leq r-1.}$}

Define $W_{a_{t_0+1}j\ll k\ll a_{t_0}j}$ by taking the summation $\sum W_{j,k}$ over all of the above integers $j,k$. In this case,
\[|\widetilde{D}(x,y)|\approx x^{A_{t_0}}\prod_{t={t_0+1}}^{r}y^{\textbf{N}\begin{bmatrix}
\begin{smallmatrix}
\cdot\\
t
\end{smallmatrix}
\end{bmatrix}}.\]
Then we have, for $|y|\leq 1$,
\begin{align*}
\int_{|x|^{a_{t_0+1}}\ll |y|\ll |x|^{a_{t_0}}}|\widetilde{D}(x,y)|^{-\frac{1}{A_r}}dx&\lesssim
\int_{|x|^{a_{t_0+1}}\ll |y|\ll |x|^{a_{t_0}}}\left(x^{A_{t_0}}\prod_{t={t_0+1}}^{r}y^{\textbf{N}\begin{bmatrix}
\begin{smallmatrix}
\cdot\\
t
\end{smallmatrix}
\end{bmatrix}}\right)^{-\frac{1}{A_r}}dx\\
&\lesssim \prod_{t=t_0+1}^ry^{-\textbf{N}\begin{bmatrix}
\begin{smallmatrix}
\cdot\\
t
\end{smallmatrix}
\end{bmatrix}/A_r}
y^{\frac{1}{a_{t_0+1}}
\left(1-\frac{A_{t_0}}{A_r}\right)}\lesssim 1
\end{align*}
since
\begin{align*}
&-\sum_{t=t_0+1}^{r}\frac{\textbf{N}\begin{bmatrix}
\begin{smallmatrix}
\cdot\\
t
\end{smallmatrix}
\end{bmatrix}}{A_r}
+\frac{1}{a_{t_0+1}}\Big(1-\frac{A_{t_0}}{A_r}\Big)
=\frac{1}{a_{t_0+1}}\cdot\frac{1}{A_r}
\Big(A_r-A_{t_0}-\sum\limits_{t=t_0+1}^{r}\textbf{N}\begin{bmatrix}
\begin{smallmatrix}
\cdot\\
t
\end{smallmatrix}
\end{bmatrix}a_{t_0+1}\Big)
\geq 0.
\end{align*}

\textrm{$\mathbf{Case~(iv)~k= a_{t_0}j + O(1)~~for~some~~1\leq t_0\leq r.}$}\\

We shall break this case into several subcases for convenience.\\

\textrm{$\mathbf{Subcase~(a)~m>0~~or~~r\geq 2}$}.

If we first integrate $|\widetilde{D}(x,y)|^{-\frac{1}{A_r}}$ with respect to $x$, then we have
\begin{align*}
&\int_{|x|\approx |y|^{\frac{1}{a_{t_0}}}}|x|^{-\frac{m}{A_r}}\prod_{l=1}^r
\prod_{r_{\nu}\in \Psi[\begin{smallmatrix}
\cdot\\
l
\end{smallmatrix}]}|y-r_\nu(x)|^{-\frac{1}{A_r}}dx\\
&\lesssim \prod_{t=t_0+1}^{r} y^{-\textbf{N}\begin{bmatrix}
\begin{smallmatrix}
\cdot\\
t
\end{smallmatrix}
\end{bmatrix}/A_r}
\int_{|x|\approx |y|^{\frac{1}{a_{t_0}}}}|x|^{-\frac{m}{A_r}}\prod_{l=1}^{t_0}
\prod_{r_{\nu}\in \Psi[\begin{smallmatrix}
\cdot\\
l
\end{smallmatrix}]}|y-r_\nu(x)|^{-\frac{1}{A_r}}dx\\
&\lesssim \prod_{t=t_0+1}^{r} y^{-\textbf{N}\begin{bmatrix}
\begin{smallmatrix}
\cdot\\
t
\end{smallmatrix}
\end{bmatrix}/A_r}
|y|^{-\frac{A_{t_0-1}}{a_{t_0}A_r}}
\int_{|x|\approx |y|^{\frac{1}{a_{t_0}}}}
\prod_{r_{\nu}\in \Psi[\begin{smallmatrix}
\cdot\\
t_0
\end{smallmatrix}]}|y-r_\nu(x)|^{-\frac{1}{A_r}}dx.
\end{align*}
By a change of variables, we have
\begin{align*}
&\int_{|x|\approx |y|^{\frac{1}{a_{t_0}}}}
\prod_{r_{\nu}\in \Psi[\begin{smallmatrix}
\cdot\\
t_0
\end{smallmatrix}]}|y-r_\nu(x)|^{-\frac{1}{A_r}}dx
\lesssim |y|^{{1}/{a_{t_0}}-{\textbf{N}\begin{bmatrix}
\begin{smallmatrix}
\cdot\\
t_0
\end{smallmatrix}
\end{bmatrix}}/{A_r}}.
\end{align*}
On the other hand, it is easy to see that
\begin{eqnarray*}
-\frac{A_{t_0-1}}{a_{t_0}A_r}-\sum_{t=t_0}^r\textbf{N}\begin{bmatrix}
\begin{smallmatrix}
\cdot\\
t
\end{smallmatrix}
\end{bmatrix}/A_r+\frac{1}{a_{t_0}}
\geq -\frac{1}{a_{t_0}}\cdot\frac{1}{A_r}\cdot A_r+\frac{1}{a_{t_0}}=0.
\end{eqnarray*}
This implies that the operator $\sum_{k= a_{t_0}j + O(1)}W_{j,k}$ is bounded on $L^1$.\\

\textrm{$\mathbf{Subcase~(b)~m=0~~and~~r=1}$}.

As classified in Section 3, all of the roots
$r_{\nu}\in
\Psi\begin{bmatrix}
\begin{smallmatrix}
\cdot\\ 1
\end{smallmatrix}
\end{bmatrix}
$
are of form $r_{\nu}(x)=C_1^{\alpha}x^{a_1}+o(x^{a_1})$. Recall that we have pointed out that two roots $r_{\nu}$ and $r_{\nu'}$ may be different even though their leading coefficients $C_1^{\alpha}$ and $C_1^{\alpha'}$ are the same.

If the leading exponent $a_1<1$, then there exist at least two different coefficients $C_1^{\alpha}$ of the leading term $x^{a_1}$. The reason is that the product $\prod_{\alpha}(y-C_1^{\alpha}x^{a_1})$ is just the polynomial associated with the edge $E$ on the boundary $\partial N(S)$ of $N(S)$ with slope $-1/a_1$. In this case, all roots $r_{\nu}$ have the same order $x^{a_1}$, and two different roots have also distance of order $x^{a_1}$. By a similar argument as in $\textbf{Subcase~(a)}$, we can prove that the integral of $|\widetilde{D}(x,y)|^{-1/A_1}$ over $|x|\approx |y|^{1/a_1}$, with respect to $x$, is bounded by a constant independent of $y$.

For $a_1>1$, by a simple change of variables, one can see that, for each root $r_{\nu}$,
\begin{equation*}
\sup_{|y|\leq 1}\int_{|x|\approx |y|^{1/a_1}} |y-r_{\nu}(x)|^{-1/a_1}dx< +\infty
\end{equation*}
with the bound independent of $y$. By H\"{o}lder's inequality, we have
\begin{eqnarray*}
& &
\sup_{|y|\leq 1}\int
\prod_{r_{\nu}\in
\Psi
\begin{bmatrix}
\begin{smallmatrix}
\cdot\\ 1
\end{smallmatrix}
\end{bmatrix}}
|y-r_{\nu}(x)|^{-1/A_1}dx
 \leq
\prod_{r_{\nu}\in
\Psi
\begin{bmatrix}
\begin{smallmatrix}
\cdot\\ 1
\end{smallmatrix}
\end{bmatrix}}
\left(
\sup_{|y|\leq 1}\int |y-r_{\nu}(x)|^{-1/a_1}dx
\right)^{a_1/A_1}<+\infty,
\end{eqnarray*}
where both integrals are taken over $|x|\approx |y|^{1/a_1}$.

Now we deal with the case $a_1=1$. By our resolution of roots in Section 3, we can decompose $W_{j,k}$ into smaller pieces of nondegenerate operators, i.e., $W_{j,k}=\sum W_{j,k,\mu_1,\cdots,\mu_N}^{\sigma_1,\cdots,\sigma_N}$. On the support of these operators, the size of $S_{xy}''$ is bounded from both above and below by positive constants.
\begin{itemize}
\item If the cluster $\Psi
\begin{bmatrix}
\begin{smallmatrix}
\alpha_1,&\cdots,&\alpha_{N-1},&\alpha_N\\
t_1,&\cdots,&t_{N-1},&t_N
\end{smallmatrix}
\end{bmatrix}$ contains at least two different roots, then $W_{j,k,\mu_1,\mu_2,\cdots,\mu_N}^{\sigma_1,\sigma_2,\cdots,\sigma_N}$ is defined as $W_{j,k}$ with insertion of the following cut-off function
\begin{equation*}
\prod_{i=1}^{N}
\Phi
\left(
\sigma_t\Big(y-C_{t_1}^{\alpha_1}x^{a_{t_1}}-
\cdots-C_{t_1,\cdots,t_i}^{\alpha_1,\cdots,\alpha_i}
x^{a_{t_1,\cdots,t_i}^{\alpha_1,\cdots,\alpha_{i-1}}}\Big)/{2^{\mu_i}}
\right).
\end{equation*}
\item Otherwise if the cluster $\Psi
\begin{bmatrix}
\begin{smallmatrix}
\alpha_1,&\cdots,&\alpha_{N-1},&\alpha_N\\
t_1,&\cdots,&t_{N-1},&t_N
\end{smallmatrix}
\end{bmatrix}$ includes only one root with its multiplicity, then $W_{j,k,\mu_1,\cdots,\mu_N}^{\sigma_1,\cdots,\sigma_N}$ is obtained by insertion of the following cut-off function into $W_{j,k}$:
\[
\prod_{i=1}^{N-1}
\Phi\left(
\sigma_t\Big(y-C_{t_1}^{\alpha_1}x^{a_{t_1}}-
\cdots-C_{t_1,\cdots,t_i}^{\alpha_1,\cdots,\alpha_i}
x^{a_{t_1,\cdots,t_i}^{\alpha_1,\cdots,\alpha_{i-1}}}\Big)/{2^{\mu_i}}
\right)\cdot \Phi\left(\sigma_N\Big(y-r_\nu(x)\Big)/{2^{\mu_N}}\right),
\]
\end{itemize}
where $r_{\nu}$ is the root (with its multiplicity) in the cluster
$\Psi
\begin{bmatrix}
\begin{smallmatrix}
\alpha_1,&\cdots,&\alpha_{N-1},&\alpha_N\\
t_1,&\cdots,&t_{N-1},&t_N
\end{smallmatrix}
\end{bmatrix}$.\\

In the first case, we must have one of the following statements.
\begin{enumerate}
\item[{\rm ($\alpha$)}] $\mu_N\gg a_{t_1,\cdots,t_{N+1}}^{\alpha_1,\cdots,\alpha_{N}}j$ with $t_{N+1}=1$ or $\mu_N\ll a_{t_1,\cdots,t_{N}, t_{N+1}}^{\alpha_1,\cdots,\alpha_{N}}$ with $t_{N+1}=\textbf{N}\begin{bmatrix}
\begin{smallmatrix}
\alpha_1,\cdots,\alpha_N\\
t_1,\cdots,t_N
\end{smallmatrix}
\end{bmatrix}.$
\item[{\rm ($\beta$)}] $a_{t_1,\cdots,t_{N+1}+1}^{\alpha_1,\cdots,\alpha_{N}}j\ll \mu_N\ll a_{t_1,\cdots,t_{N+1}}^{\alpha_1,\cdots,\alpha_{N}}j$
   ~for some $1\leq t_{N+1} \leq \textbf{N}\begin{bmatrix}
   \begin{smallmatrix}
\alpha_1,\cdots,\alpha_N\\
t_1,\cdots,t_N
\end{smallmatrix}
\end{bmatrix}-1$.
\item[{\rm ($\gamma$)}] $\mu_N= a_{t_1,\cdots,t_{N+1}}^{\alpha_1,\cdots,\alpha_{N}}j+O(1)$ and on the support of $W_{j,k,\mu_1,\cdots,\mu_N}^{\sigma_1,\cdots,\sigma_N}$ it is true that for all $\alpha$
   \begin{equation*}
    \Big|
    y-\sum_{i=1}^{N} C_{t_1,\cdots,t_i}^{\alpha_1,\cdots,\alpha_i}
     x^{a_{t_1,\cdots,t_i}^{\alpha_1,\cdots,\alpha_{i-1}}}
     -
     C_{t_1,\cdots,t_N,t_{N+1}}^{\alpha_1,\cdots,\alpha_N,\alpha}
     x^{a_{t_1,\cdots,t_N,t_{N+1}}^{\alpha_1,\cdots,\alpha_N}}
    \Big|
    \approx
    2^{\mu_N}.
   \end{equation*}
\end{enumerate}
Note that $a_{t_1}=t_1=1$ here and we do not write it explicitly to lighten the notation. \\

Assume we are in the case ($\beta$). On the support of $W_{j,k,\mu_1,\mu_2,\cdots,\mu_N}^{\sigma_1,
\sigma_2,\cdots,\sigma_N}$,
the size of $\widetilde{D}$ (equal to $D$) is equivalent to
\begin{align*}
&
\prod_{\alpha\neq \alpha_1}
2^{k\textbf{N}\begin{bmatrix}
\begin{smallmatrix}
\alpha \\ 1
\end{smallmatrix}
\end{bmatrix}}
\prod_{t<t_2}2^{j\textbf{N}\begin{bmatrix}
\begin{smallmatrix}
\alpha_1 & \cdot\\
t_1 & t
\end{smallmatrix}
\end{bmatrix}a_{t_1,t}^{\alpha_1}}
\prod_{t>t_2}2^{\mu_1\textbf{N}\begin{bmatrix}
\begin{smallmatrix}
\alpha_1 & \cdot\\
t_1 & t
\end{smallmatrix}
\end{bmatrix}}
\prod_{\alpha\neq \alpha_2}2^{\mu_1\textbf{N}\begin{bmatrix}
\begin{smallmatrix}
\alpha_1 & \alpha\\
t_1 & t_2
\end{smallmatrix}
\end{bmatrix}}
\cdots\\
&~~~~~~~\prod_{t\leq t_{N+1}}2^{j\textbf{N}\begin{bmatrix}
\begin{smallmatrix}
\alpha_1 & \alpha_2 & \cdots & \alpha_N &\cdot\\
t_1 & t_2 & \cdot & t_N & t
\end{smallmatrix}
\end{bmatrix}a_{t_1,t_2,\cdots, t_N, t}^{\alpha_1, \alpha_2, \cdots, \alpha_N}}
\prod_{t>t_{N+1}}2^{\mu_N\textbf{N}\begin{bmatrix}
\begin{smallmatrix}
\alpha_1 & \alpha_2 & \cdots & \alpha_N &\cdot\\
t_1 & t_2 & \cdot & t_N & t
\end{smallmatrix}
\end{bmatrix}}.
\end{align*}
On the other hand, by our resolution of roots, one has
\begin{align}\label{Section 4 relation mu1 mu2 muN}
&2^k\approx 2^{ja_{t_1}}, \quad 2^{\mu_1}\approx 2^{ja_{t_1,t_2}^{\alpha_1}}, \cdots, 2^{\mu_{N-1}}\approx 2^{ja_{t_1,t_2,\cdots,t_{N}}^{\alpha_1,\alpha_2,\cdots,\alpha_{N-1}}}.
\end{align}
Hence the size of $\widetilde{D}$ is not less than a constant multiple of
\begin{equation*}
2^{ja_{t_1,t_2,\cdots,t_{N+1}}^{\alpha_1,\alpha_2,\cdots,\alpha_{N}}
\Big(
\textbf{N}\begin{bmatrix}
\begin{smallmatrix}
\cdot \\ 1
\end{smallmatrix}
\end{bmatrix}
-\sum_{t>t_N}
\textbf{N}\begin{bmatrix}
\begin{smallmatrix}
\alpha_1,&\cdots,&\alpha_N,&\cdot \\
t_1,&\cdots,&t_N,t
\end{smallmatrix}
\end{bmatrix}
\Big)
2^{\mu_N
\sum_{t>t_N}
\textbf{N}\begin{bmatrix}
\begin{smallmatrix}
\alpha_1,&\cdots,&\alpha_N,&\cdot \\
t_1,&\cdots,&t_N,&t
\end{smallmatrix}
\end{bmatrix}
}
}.
\end{equation*}
By the size estimate, the horizontal cross sections of the support of $W_{j,k,\mu_1,\mu_2,\cdots,\mu_N}^{\sigma_1,\sigma_2,
\cdots,\sigma_N}$ have size $\lesssim 2^{\mu_N}$. Hence
\begin{align}\label{section 4 boundedness Wjk on L1}
\sum_{\mu_N}\int_{\mathbb{R}}
\left|W_{j,k,\mu_1,\mu_2,\cdots,\mu_N}^{\sigma_1,\sigma_2,
\cdots,\sigma_N}f(x)\right|dx
&\lesssim  \sum_{\mu_N}\int_{\mathbb{R}}
2^{\mu_{N}}\sup\left(|\widetilde{D}(x,y)|^{-\frac{1}{A_1}}\right)
|f(y)|dy
\lesssim  \int_{\mathbb{R}}|f(y)|dy,
\end{align}
where the above summations are taken over $\mu_N$ in the range $(\beta)$. Note
that the number of $\mu_1,\cdots,\mu_{N-1}$ is bounded by a constant independent of $j,k$. Therefore $W_{j,k}$ is bounded on $L^1$.

Now we assume the cluster $\Psi
\begin{bmatrix}
\begin{smallmatrix}
\alpha_1,&\cdots,&\alpha_{N-1},&\alpha_N\\
t_1,&\cdots,&l_{N-1},&t_N
\end{smallmatrix}
\end{bmatrix}$
contains only one root with certain multiplicity. It suffices to study $\mu_N$ in the range $\mu_N\ll a_{t_1,\cdots,t_{N+1}}^{\alpha_1,\cdots,\alpha_N}j$. In fact, by our assumption, there is only one leading exponent $a_{t_1,\cdots,t_{N+1}}^{\alpha_1,\cdots,\alpha_N}$. Hence the case $\mu_N\geq a_{t_1,\cdots,t_{N+1}}^{\alpha_1,\cdots,\alpha_N}j-K$ is essentially the same as $(\alpha)$ discussed above. Since $D$ is not of the form $(y-r_{\nu}(x))^d$ for some root $r_{\nu}(x)=C_1x+o(x)$, we can repeat the above argument step by step to prove that $W_{j,k}$ is bounded on $L^1$.

Combining all above results, we have completed the proof of $L^1$ boundedness for $W_{j,k,\mu_1,\mu_2,\cdots,\mu_N}^{\sigma_1,\sigma_2,\cdots,\sigma_N}$ when the damping function $D$ is not of the form $(y-r_\nu(x))^d$.\\

It remains to consider the case $D(x,y)=(y-r_\nu(x))^d$ with
$r_\nu(x)=C_1x+o(x)$. Here $C_1\neq 0$ is a real number. By dilation, we may assume $C_1=1$. Each root $r_{\nu}$ in $\Psi\begin{bmatrix}
\begin{smallmatrix}
\cdot\\
1
\end{smallmatrix}
\end{bmatrix}$ is a real Puiseux series. In a small neighborhood of the origin, $r_{\nu}'(x)$ is bounded from above and below by fixed positive constants. Hence the inverse $r_{\nu}^{-1}$ is well defined. Now we are going to prove that
\begin{equation}\label{section 4 L1H1 boundedness}
\|W_{j,k}f\|_{L^1}\leq C\|f\|_{H_E^1(I_k)}, \quad \Re(z)=-\frac{1}{A_1}
\end{equation}
for each pair $(j,k)\in\mathbb{Z}^2$ satisfying $j-N_0\leq k\leq j+N_0$. Here $N_0$ denotes a large positive integer.

To prove the above $H^1_E\rightarrow L^1$ estimate, we first need the following $L^2$ estimate
\begin{equation*}
\|W_{j,k,}f\|_{L^2}\leq C\|f\|_{L^2}, \quad j-N_0\leq k\leq j+N_0, \quad\Re(z)=-\frac{1}{A_1}.
\end{equation*}
Define $W_{j,k,l}^{\sigma}$ by inserting $\Phi\left(\sigma\frac{y-r_\nu(x)}{2^l}\right)$ into $W_{j,k}$. Then
$W_{j,k}=\sum_l\sum_{\sigma=\pm}W_{j,k,l}^{\sigma}.$
By the Schur test, it suffices to consider those $l\in\mathbb{Z}$ such that
\begin{equation}\label{section 4 large l}
2^{l\textbf{N}\begin{bmatrix}
\begin{smallmatrix}
\cdot\\
1
\end{smallmatrix}
\end{bmatrix}}\geq\left(|\lambda|2^{kB_1}\right)^{-\frac{\textbf{N}\begin{bmatrix}
\begin{smallmatrix}
\cdot\\
1
\end{smallmatrix}
\end{bmatrix}}{\textbf{N}\begin{bmatrix}
\begin{smallmatrix}
\cdot\\
1
\end{smallmatrix}
\end{bmatrix}+2}}.
\end{equation}
By Lemma \ref{Lemma operator vers van der Corput}, we have
\[\left\|W_{j,k,l}^{\sigma}\right\|\leq C\left(|\lambda|2^{l\textbf{N}\begin{bmatrix}
\begin{smallmatrix}
\cdot\\
1
\end{smallmatrix}
\end{bmatrix}}2^{kB_1}\right)^{-\frac{1}{2}}2^{-l}.\]
Hence we have
\begin{align*}
\Big\|\sum_{l}W_{j,k,l}^{\sigma}\Big\|&\leq C\sum_{l}\left(|\lambda|2^{l\textbf{N}\begin{bmatrix}
\begin{smallmatrix}
\cdot\\
1
\end{smallmatrix}
\end{bmatrix}}2^{kB_1}\right)^{-\frac{1}{2}}2^{-l}\leq C,
\end{align*}
where the summation is taken over all $l$ satisfying
(\ref{section 4 large l}). Hence the above $L^2$ estimate is true.

In what follows, it will be shown that $\|W_{j,k}a\|_{L^1}\leq C\|a\|_{H_E^1(I_k)}$ for all $H_E^1(I_k)$ atoms $a$. Assume $a$ is an $H_E^1(I_k)$ atom. Then there exists an interval $I\subseteq I_k$ for which we have
(i) $\supp(a)\subset I$; (ii) $\|a\|_{L^\infty}\leq |I|^{-1}$;
(iii $\int_{I}e^{i\lambda S(r_{\nu}^{-1}(c_I),y)}a(y)dy=0$.

Write $I=(c_I-\delta, c_I+\delta)$ for some $\delta>0$. Now we define an interval $J$ associated with $I$. If $|I|\geq \left(|\lambda|2^{kB_1}\right)^{{-\frac{1}{\textbf{N}\begin{bmatrix}
\begin{smallmatrix}
\cdot\\
1
\end{smallmatrix}
\end{bmatrix}+2}}}$, then set $J=(r_{\nu}^{-1}(c_I)-2\delta, r_{\nu}^{-1}(c_I)+2\delta)$. Otherwise we define $J$ as
\[J=\Big(r_{\nu}^{-1}(c_I)-2\left(|\lambda|2^{kB_1}\right)^{
{-\frac{1}{\textbf{N}\begin{bmatrix}
\begin{smallmatrix}
\cdot\\
1
\end{smallmatrix}
\end{bmatrix}+2}}}, r_{\nu}^{-1}(c_I)+2\left(|\lambda|2^{kB_1}\right)^{{-\frac{1}{\textbf{N}\begin{bmatrix}
\begin{smallmatrix}
\cdot\\
1
\end{smallmatrix}
\end{bmatrix}+2}}}\Big).\]
For the interval $J$, we claim that there exists a constant $C>0$ such that
\[\|W_{j,k}a\|_{L^1(J)}\leq C\|a\|_{H_E^1(I_k)}\]
for all $H_E^1(I_k)$ atoms $a$. In fact, by Cauchy-Schwarz's inequality, we have
\begin{align*}
\|W_{j,k}a\|_{L^1(r_{\nu}^{-1}(c_I)-2\delta, r_{\nu}^{-1}(c_I)+2\delta)}&\leq (4\delta)^{\frac{1}{2}}\|W_{j,k}a\|_{L^2}
\leq C\delta^{\frac{1}{2}}\|a\|_{L^2}
\leq C.
\end{align*}
On the other hand, we also claim that there exists a constant $C$ such that
\[\|W_{j,k}a\|_{L^1(r_{\nu}^{-1}(c_I)-2\widetilde{\delta}, r_{\nu}^{-1}(c_I)+2\widetilde{\delta})}\leq C,\]
with $\widetilde{\delta}=\left(|\lambda|2^{kB_1}\right)^{{-\frac{1}{\textbf{N}\begin{bmatrix}
\begin{smallmatrix}
\cdot\\
1
\end{smallmatrix}
\end{bmatrix}+2}}}$. Consider the integral kernel
$
K(x,y)=\left|D(x,y)\right|^{-\frac{1}{A_1}}
\chi_{J}(x)
\chi_I(y).
$
It is easy to see that
\begin{align*}
&\sup_x\int_{\mathbb{R}}|K(x,y)|dy\leq C\left(|\lambda|2^{kB_1}\right)^{\frac{1}{\textbf{N}\begin{bmatrix}
\begin{smallmatrix}
\cdot\\
1
\end{smallmatrix}
\end{bmatrix}+2}}|I|,\\
&\sup_y\int_{\mathbb{R}}|K(x,y)|dx\leq C\left(|\lambda|2^{kB_1}\right)^{\frac{1}{\textbf{N}\begin{bmatrix}
\begin{smallmatrix}
\cdot\\
1
\end{smallmatrix}
\end{bmatrix}+2}}\left(|\lambda|2^{kB_1}\right)^{-\frac{1}{\textbf{N}\begin{bmatrix}
\begin{smallmatrix}
\cdot\\
1
\end{smallmatrix}
\end{bmatrix}+2}}\leq C.
\end{align*}
By the Schur test, one has
\begin{align*}
\|W_{j,k}a\|_{L^1(r_{\nu}^{-1}(c_I)-2\widetilde{\delta}, r_{\nu}^{-1}(c_I)+2\widetilde{\delta})}&\leq C\sqrt{\widetilde{\delta}}\|W_{j,k}a\|_{L^2(r_{\nu}^{-1}(c_I)-2\widetilde{\delta}, r_{\nu}^{-1}(c_I)+2\widetilde{\delta})}\\
&\leq C\sqrt{\widetilde{\delta}}\frac{1}{\sqrt{\widetilde{\delta}}}|I|^{\frac{1}{2}}|I|^{-\frac{1}{2}}\\
&\leq C.
\end{align*}
Thus our above claim is true. Now it remains to show that
$\|W_{j,k}a\|_{L^1(J^c)}\leq C.$
We need the following inequality.
\[\sup_{y\in I}\int_{J^c}\left|\Phi\left(\frac{x}{2^j}\right)\Phi\left(\frac{y}{2^k}\right)|\widetilde{D}(x,y)|^z\varphi(x,y)-
\Phi\left(\frac{x}{2^j}\right)\Phi\left(\frac{c_I}{2^k}\right)|\widetilde{D}(x,c_I)|^z\varphi(x,c_I)\right|dx\leq C.\]
Define $M_{j,k}(x,y)
=\Phi
\left(\frac{x}{2^j}\right)
\Phi\left(\frac{y}{2^k}\right)\varphi(x,y)$. By the mean value theorem, it follows immediately from $k\geq 0$ that
\begin{equation*}
|M_{j,k}(x,y)-M_{j,k}(x,c_I)|
\leq
\|\Phi\|_{\infty}\|\Phi'\|_{\infty}\frac{|I|}{2^k}\|\varphi\|_{\infty}+
\|\Phi\|_{\infty}^2\|\partial_y\varphi\|_{\infty}|I|\leq C|I|/2^k.
\end{equation*}
Observe that for $x\in J^c$, $\Phi(x/2^j)\neq 0$, $ y\in I$, it is true that
$|y-r_\nu(x)|\approx |r_{\nu}^{-1}(y)-x|\geq C|I|.$
Also, for $y\in I$ and $x\in J^c$,
$\left||\widetilde{D}(x,y)|^z-|\widetilde{D}(x,c_I)|^z\right|\leq C|I||x-r_{\nu}^{-1}(c_I)|^{-2}.$
It is clear that for $x\in J^c$ we have
$\left||\widetilde{D}(x,y)|^z\right|\leq \left|x-r_{\nu}^{-1}(c_I)\right|^{-1}$.
Then we have
\begin{align*}
&\int_{J^c\cap I_j}\left|M_{j,k}(x,y)|\widetilde{D}(x,y)|^z-M_{j,k}(x,c_I)|\widetilde{D}(x,c_I)|^z\right|dx\\
&\leq C\int_{|I|\leq |x-r_{\nu}^{-1}(c_I)|\leq C2^j}\frac{|I|}{2^k}|x-r_{\nu}^{-1}(c_I)|^{-1}dx+\int_{|x-r_{\nu}^{-1}(c_I)|\geq |I|}\frac{|I|}{|x-r_{\nu}^{-1}(c_I)|^2}dx\leq C.
\end{align*}
By the cancellation property of $a$, we have
\begin{align*}
\left|\int_{I}e^{i\lambda S(x,y)}a(y)dy\right|&=\left|\int_{I}\left[e^{i\lambda (S(x,y)-S(x,c_I))}-e^{i\lambda (S(r_{\nu}^{-1}(c_I),y)-S(r_{\nu}^{-1}(c_I),c_I))}\right]a(y)dy\right|\\
&\leq C|\lambda||I||x-r_{\nu}^{-1}(c_I)|^{\textbf{N}\begin{bmatrix}
\begin{smallmatrix}
\cdot\\
1
\end{smallmatrix}
\end{bmatrix}+1}2^{kB_1}, \quad x\in J^c\cap I_j.
\end{align*}
Hence
\begin{eqnarray*}
& &\int_{x\in J^c:~|x-r_{\nu}^{-1}(c_I)|\leq \eta}\left|\int_{I}e^{i\lambda S(x,y)}a(y)dy\right|
\left| M_{j,k}(x,c_I)|\widetilde{D}(x,c_I)|^z \right| dx\\
&\leq&
C\int_{x\in J^c:~|x-r_{\nu}^{-1}(c_I)|\leq \eta}
|x-r_{\nu}^{-1}(c_I)|^{-1}
\left|\int_{I}e^{i\lambda S(x,y)}a(y)dy\right|dx\\
&\leq&
C\int_{x\in J^c:~|x-r_{\nu}^{-1}(c_I)|\leq \eta}
|\lambda||I|\Big|x-r_{\nu}^{-1}(c_I)\Big|^
{\textbf{N}
\begin{bmatrix}
\begin{smallmatrix}
\cdot\\ 1
\end{smallmatrix}
\end{bmatrix}
}
2^{kB_1}dx\\
&\leq &
C|\lambda||I|\eta^{\textbf{N}
\begin{bmatrix}
\begin{smallmatrix}
\cdot\\ 1
\end{smallmatrix}
\end{bmatrix}
+1}
2^{kB_1}
\end{eqnarray*}
Choose $\eta>0$ such that $|\lambda||I|\eta^{\textbf{N}\begin{bmatrix}
\begin{smallmatrix}
\cdot\\
1
\end{smallmatrix}
\end{bmatrix}+1}2^{kB_1}=1.$
Let $\mu=\max\{\eta,|I|\}=\max\{\eta,2\delta\}$. Then we have proved $\|W_{j,k}a\|_{L^1(|x-r_{\nu}^{-1}(c_I)|\leq \mu)}\leq C$. For $x\in [2^{j-1},2^{j+1}]$ and
$|x-r_{\nu}^{-1}(c_I)|>\mu$, it follows from Lemma \ref{Lemma uniform L2 estimate} that
\begin{align*}
&\sum_{\mu\leq 2^l \lesssim 2^j}\int_{2^{l-1}\leq|x-r_{\nu}^{-1}(c_I)|\leq 2^{l+1}}|x-r_{\nu}^{-1}(c_I)|^{-1}\left|\int_{I}e^{i\lambda S(x,y)}a(y)dy\right|dx\\
&\leq C \sum_{\mu\leq 2^l \lesssim 2^j} \int_{\frac{1}{2}\leq |x|\leq 2}\Big|\int_{[-1, 1]}e^{i\lambda S(2^lx+r_{\nu}^{-1}(c_I), \delta y+c_I)}\delta a(\delta y+c_I)dy\Big|dx\\
&\leq C \sum_{\mu\leq 2^l \lesssim 2^j}
\Big(
|\lambda|\delta 2^{l
\Big(\textbf{N}\begin{bmatrix}
\begin{smallmatrix}
\cdot\\1
\end{smallmatrix}
\end{bmatrix}+1\Big)}
2^{kB_1}\Big)^{-1/2}<+\infty.
\end{align*}
Combining all above results, we obtain the $H_E^1\to L^1$ boundedness for $W_{j,k}$. Hence the proof of the theorem is complete.
\end{proof}

\section{Proof of Theorem \ref{section 1 Main theorem}}
In this section, we shall prove Theorem 1 by interpolating the damping estimates in Sections 3 and 4.

\begin{proof}
By a duality argument, we can assume $A_r\geq B_r$ in Theorem \ref{section 3 Thm main result}. Hence it suffices to establish the $L^p$ estimate (\ref{sec-main-est}) for $A_r\geq B_r$.

We first consider the case when the damping factor $D$ in (\ref{section 3 Def damping D}) is not of the form (\ref{section 4 translation invar damp fac}).
As in Section 4, it is more convenient to consider the following four cases separately. \\

\textbf{Case (i) $k\gg a_1j$}

In this case, we define $\widetilde{D}(x,y)=D(x,y)$ and $W_z^{(3)}=\sum_{j,k}W_{j,k}$ with the summation taken over all above $j,k$. By Theorems \ref{Thm L2 damping decay estiamte} and \ref{Thm Main theorem Sec 4}, we have the $L^2$ decay estimate
\begin{eqnarray}\label{section 5 L2 damping estimate}
\|W_z^{(3)}f\|_{L^2}
&\leq&
C(1+|z|^2)|\lambda|^{-\frac{1}{2(1+B_r)}} \|f\|_{L^2},
~~~\Re(z)=\frac{A_r-B_r}{2A_r}\cdot\frac{1}{1+B_r}
\end{eqnarray}
and the $L^1\rightarrow L^1$ estimate
\begin{eqnarray}\label{section 5 L1 damping estimate}
\|W_z^{(3)}f\|_{L^1}
&\leq&
C\|f\|_{L^1},~~~\Re(z)=-\frac{1}{A_r},
\end{eqnarray}
where the above constants $C$ are independent of $\lambda$ and $f$. We shall point out that the $L^2$ decay estimate in Theorem \ref{Thm L2 damping decay estiamte}
is also true for the three different parts $W_z^{(1)}$, $W_z^{(2)}$ and $W_z^{(3)}$, as defined in Section 4. In the proof of Theorem \ref{Thm L2 damping decay estiamte}, we have decomposed $W_z$ suitably into different parts and then proved $L^2$ desired estimate for each of those parts. Hence the above $L^2$ estimate (\ref{section 5 L2 damping estimate}) for $W_z^{(3)}$ also holds.

By Stein's complex interpolation (see Stein \cite{stein}), we see that
\begin{equation}\label{section 5 desired Lp estimate}
\|W_z^{(3)}f\|_{L^{p_r}}
\leq
C
|\lambda|^{-\frac{1}{A_r+B_r+2}}
\|f\|_{L^{p_r}},~~~p_r=\frac{A_r+B_r+2}{A_r+1},
~~~\Re(z)=0.
\end{equation}

\textbf{Case (ii) $k\ll a_rj$}

In this case, we shall modify our definition $\widetilde{D}$ in Section 4. Indeed, let $\widetilde{D}(x,y)=|x|^{A_r}$ and $W_z^{(1)}=\sum W_{j,k}$ with the summation taken over all $j,k$ satisfying $k\ll a_rj$. It is easy to see that $|D(x,y)|\approx |x|^{A_r}$ in the support of $W_z^{(1)}$. With the same argument as in Section 3, one can show that the $L^2$ decay estimate (\ref{section 5 L2 damping estimate}) is also true for $W_z^{(1)}$. But $W_z^{(1)}$ does \textbf{not} satisfy the $L^1\rightarrow L^1$ estimate
(\ref{section 5 L1 damping estimate}). However, we can replace (\ref{section 5 L1 damping estimate}) by the following $L^1\rightarrow L^{1,\infty}$ estimate:
\begin{eqnarray}
\|W_z^{(1)}f\|_{L^{1,\infty}}
&\leq&
C\|f\|_{L^1},~~~\Re(z)=-\frac{1}{A_r}.
\end{eqnarray}
Using the interpolation technique in Lemma \ref{section 2 stein-weiss interpolation}, we obtain the desired estimate (\ref{section 5 desired Lp estimate}) for $W_z^{(1)}$.\\

\textbf{Case (iii) $a_{t+1}j\ll k\ll a_tj$ for some $1\leq t \leq r-1$.}

We define $W_z^{(2)}=\sum W_{j,k}$ with the summation taken over all above $j,k$. In this case, $\widetilde{D}(x,y)=D(x,y)$. As shown in Sections 3 and 4, we see that the same estimates as in Case (i) are also true for $W_z^{(2)}$. An interpolation yields the desired $L^p$ estimate.\\

\textbf{Case (iv) $k= a_tj+O(1)$ for some $1\leq t \leq r$.}

Since we have assumed that $D$ is not of the form (\ref{section 4 translation invar damp fac}), all estimates in Case (i) are true in this case. The desired $L^p$ decay estimate follows immediately.\\

Now we turn to the case in which $D(x,y)=(y-r_\nu(x))^d$ for some positive integer $d$ and a nontrivial root $r_\nu(x)=C_1x+o(x)$. Here $C_1$ is a nonzero real number. In this case, we need only take care of \textbf{Case~(iv)} as above since other cases can be treated in the same way.

Now we assume $|j-k|\leq N_0$ for some large $N_0>0$. Let the damping factor $\widetilde{D}$ be defined as in (\ref{section 4 def damping factor}). Similarly, we use $W_z^{(2)}$ to denote the summation $\sum_{|j-k|\leq N_0}W_{j,k}$. Then $W_z^{(2)}$ still satisfies (\ref{section 5 L2 damping estimate}). In Theorem \ref{Thm Main theorem Sec 4}, we have the following appropriate substitute of the $L^1$ estimate $(\ref{section 5 L1 damping estimate})$:
\begin{equation}
\|W_z^{(2)}f\|_{L^1}
\leq
C \|f\|_{H_E^1(I_k)},~~~\Re(z)=-\frac{1}{N[\begin{smallmatrix}
\cdot\\
1
\end{smallmatrix}]}.
\end{equation}
It should be pointed out that $A_1=N[\begin{smallmatrix}
\cdot\\
1
\end{smallmatrix}]$ in this case. Using the Fefferman-Stein sharp function in \cite{feffermanstein} (see also \cite{PS1986} for a variant of this function), we see that an interpolation between $L^2\rightarrow L^2$ and $H_E^1\rightarrow L^1$  yields the desired estimate (\ref{section 5 desired Lp estimate}).\\

Combining all above results, we have completed the proof of Theorem \ref{section 3 Thm main result}. By interpolation, one can see that the decay estimate (\ref{section 1 main estimate}) is sharp only if $(k,l)$ is a vertex of $N(S)$. Hence Theorem \ref{section 1 Main theorem} follows from Theorem \ref{section 3 Thm main result}.

\end{proof}

\noindent{\bf Acknowledgements.} This work was supported in part by the National Natural Science Foundation of China under Grant No. 11701573. We would like to express our gratitude to Professor Xiaochun Li for his valuable comments and warm encouragement.

\end{document}